\documentclass[reqno]{amsart}
\usepackage{yhmath,amsmath,amsfonts,amssymb,enumerate,amsthm,appendix,caption,lipsum,}
\usepackage{mathrsfs}
\usepackage{fancyhdr}
\usepackage[english]{babel}
\usepackage{cmap,mathtools}
\usepackage{graphics} %% add this and next lines if pictures should be in esp format
\usepackage{epsfig} %For pictures: screened artwork should be set up with an 85 or 100 line screen
\usepackage{graphicx}\usepackage{epstopdf}%This is to transfer .eps figure to .pdf figure; please compile your paper using PDFLeTex or PDFTeXify.
\usepackage[pagewise]{lineno}%\linenumbers
\usepackage[margin=0.85in]{geometry}
\usepackage{cite}
\usepackage[utf8]{inputenc}
\usepackage{xcolor}
\usepackage{hyperref}
\usepackage[hyperpageref]{backref}
\hypersetup{
    colorlinks=true,
    linkcolor=myblue,
    filecolor=magenta,      
    urlcolor=mygreen,
    citecolor=mygreen,
}
\definecolor{mygreen}{rgb}{0.01,0.6,0.2}
\definecolor{myblue}{rgb}{0.01, 0.18, 1.0}
\usepackage[foot]{amsaddr}%email at bottom first page
\newtheorem{theorem}{Theorem}
\newtheorem{proposition}[theorem]{Proposition}
\newtheorem{lemma}[theorem]{Lemma}

\theoremstyle{definition}
\newtheorem{definition}[theorem]{Definition}
\newtheorem{remark}[theorem]{Remark}
\newtheorem{example}[theorem]{Example}

\numberwithin{equation}{section}
\numberwithin{theorem}{section}
\numberwithin{equation}{section}
\numberwithin{theorem}{section}
\title[Global branching for semilinear fractional Laplace]{Global branching for semilinear fractional Laplace with sublinear nonlinearity}

\author{Jefferson Abrantes}
\address[J. Abrantes]{Universidade Federal de Campina Grande, Unidade Acad\^{e}mica de Matem\'{a}tica,
CEP: 58429-970, Campina Grande - PB, Brazil.}
\email{jefferson@mat.ufcg.edu.br}

\usepackage{orcidlink}
\author{Rohit Kumar\,\orcidlink{0009-0001-6494-6407}}
\address[R. Kumar]{Tata Institute of Fundamental Research-Centre For Applicable Mathematics, Bangalore 560065, India.}
\email{rohit1.iitj@gmail.com, rohit24@tifrbng.res.in}

\author{Abhishek Sarkar\,\orcidlink{0000-0002-9097-7784}}
 \address[A. Sarkar]{Department of Mathematics, Indian Institute of Technology Jodhpur, NH 65, Jodhpur, Rajasthan - 342030, India.}
 \email{\tt  abhisheks@iitj.ac.in}

\keywords{Fractional Laplace operator, sublinear nonlinearity, variational methods, Classical Linking Theorem, multiplicity of solutions, positive solutions, regularity of solutions. \\
\phantom{aa} 2010 AMS Subject Classification: 35A15, 35B38, 35B09, 35D30, 58E05, 35J61, 35R15. }
\usepackage[T1]{fontenc}
\usepackage{enumitem}
\usepackage{mathtools,bm}
\usepackage{cmap}
\usepackage[hyperpageref]{backref}
\usepackage{bigints}
\usepackage{esint}

\DeclarePairedDelimiter\abs{\lvert}{\rvert}
\DeclarePairedDelimiter\norm{\lVert}{\rVert}
\makeatletter
\let\oldnorm\norm
\def\norm{\@ifstar{\oldnorm}{\oldnorm*}}

\makeatletter

\newcommand{\sig} {\sigma}

\DeclareMathAlphabet{\mathpzc}{T1}{pzc}{m}{it}

\def\dr{{\rm d}r}
\definecolor{tealgreen}{rgb}{0.0, 0.51, 0.5}

\definecolor{green(html/cssgreen)}{rgb}{0.0, 0.5, 0.0}

\def\Ds2{{{\mathcal D}^{s,2}(\R^N)}}
\def\R{{\mathbb R}}
\def\dx{{\rm d}x}
\def\dy{{\rm d}y}
\def\dxy{{\rm d}x {\rm d}y}
\def\ds{{\rm d}s}
\def\dz{{\rm d}z}

\def\Gag2{\iint_{\R^{2N}} \frac{|u(x)-u(y)|^2}{|x-y|^{N+2s}}\ \dxy}
\begin{document}

\begin{abstract}
This article investigates the existence, nonexistence, and multiplicity of positive solutions to the sublinear fractional elliptic problem \eqref{main-problem}. We begin by establishing several a priori estimates that provide regularity results and describe the qualitative behavior of solutions. A critical threshold level for the parameter $\lambda$ is identified, which plays a crucial role in determining the existence or nonexistence of solutions. The sub- and supersolution method is employed to obtain a weak solution. Furthermore, we establish a relation between the local minimizers of $\mathcal{D}^{s,2}(\R^N)$ versus $C(\mathbb{R}^N; 1+|x|^{N-2s})$. Combining these results with the Classical Linking Theorem, we demonstrate the existence of at least two distinct positive weak solutions to \eqref{main-problem}. This work extends the results of Yang, Abrantes, Ubilla, and Zhou (J. Differential Equations, 416:159–189, 2025) to the nonlocal setting, i.e., when $s \in (0,1)$. Several technical challenges arise in this framework, such as the lack of a standard comparison principle in $\R^N$ in the fractional setting.
\end{abstract}

\maketitle
\section{Introduction}
In recent years, there has been growing interest in nonlocal diffusion equations, especially associated with the fractional Laplacian operator. This operator arises naturally in a wide range of applications such as image processing (see Gilboa and Osher~\cite{Gilboa-Osher-2008}), crystal dislocation and materials science (see Dipierro, Palatucci and Valdinoci~\cite{Dipierro-Palatucci-Valdinoci-CMP-2015}), as well as in models from mathematical biology (see Bueno and Grau~\cite{Bueno-2014}) and probability theory, where it provides a natural description of symmetric L\'{e}vy jump processes (see Applebaum~\cite{Applebaum-2009-book}, Bass~\cite{Bass-book-1995}). The study of the fractional Laplace operator gained much attention after the pioneering work of Caffarelli and Silvestre~\cite{Cafarelli-Silvestre-ExT-2007}. For an introductory understanding of fractional Laplace operators, we refer to the concise exposition in~\cite{Bisci-Redulesco-Servadei-2016-Book, Hitchiker-guide-2012, Abatangelo-Valdinoci-Getting-aquainted,Chen-Li-Fractional-Laplacian}. For additional problems involving nonlocal operators, see~\cite{Quaas-2025-CVPDE, Garofalo-fractional-thoughts, Abatangelo-Dipierro-Valdinoci-fractional-world-book, Musina-Part-1,Musina-Part-2,Musina-part-3, Biswas-Kumar-2024-Semipositone, HKMS-Asymptotic-Analysis, Kumar-Sarkar-Mukherjee-DIE, KMS-IJPAM, Kumar-Sarkar-TMNA, Kumar-Ortega-2025}.       

Motivated by these applications and by the growing interest in fractional elliptic equations, in this paper we investigate the \emph{existence}, \emph{nonexistence}, and \emph{multiplicity} results for positive solutions of problems with sublinear nonlinearities and nonnegative weight functions. More precisely, we consider the following problem:
\begin{align}\label{main-problem}
    \begin{cases}
(-\Delta)^s u = \lambda h(x)f(u) & \text{in } \mathbb{R}^N ,\\
u > 0 & \text{in } \mathbb{R}^N, \\
u \in \mathcal{D}^{s,2}(\mathbb{R}^N),
\end{cases} \tag{$P_\lambda^s$}
\end{align}
where $s\in (0,1), N>2s$, $\lambda > 0$ is a parameter. The nonlinearity $f$ satisfies the following assumptions:
\begin{enumerate}[label={($\bf f{\arabic*}$)}]
    \setcounter{enumi}{0}
    \item \label{f1}  $f: \mathbb{R} \to \mathbb{R}$ is continuous, $f(t) = 0$ for $t \leq 0$ and $f(t) > 0$ for $t > 0$. \vspace{0.2 cm}
    \item \label{f2} $\displaystyle \lim_{t \to \infty} \frac{f(t)}{t} = 0$, i.e. $f$ is sublinear at $\infty$.\vspace{0.2 cm}
    \item \label{f3} $f$ is non-decreasing in $\mathbb{R}_+=(0,\infty)$ and locally Lipschitz in $\mathbb{R}$.\vspace{0.2 cm}
    \item \label{f4} From \ref{f1}, \ref{f2} and \ref{f3}, there exists a constant $C_0 > 0$ such that 
    \begin{align}\label{eq-f4}
        f(t) \leq C_0 t, \quad \forall t \geq 0.
    \end{align}
   \end{enumerate}
Regarding the nonnegative function \( h \), we assume:
\begin{enumerate}[label={($\bf H$)}]
    \setcounter{enumi}{0}
    \item \label{h} The function $h \in L^\infty(\mathbb{R}^N) \cap L^1(\mathbb{R}^N)$.
\end{enumerate}
Moreover, there exists a continuous function $P \in C(\mathbb{R}_+, \mathbb{R}_+)$ such that  
    \begin{align}\label{h-bound}
        0 < h(x) \leq P(|x|), \quad \forall x \in \mathbb{R}^N \setminus \{0\},
    \end{align}
and satisfying:
\begin{enumerate}[label={($\bf P1$)}]
    \setcounter{enumi}{0}
    \item \label{P1} $ \int_{\R^N} P(|y|)|x-y|^{2s-N}\,\dy \in L^\infty(\R^N)$.
\end{enumerate}
\begin{enumerate}[label={($\bf P2$)}]
    \setcounter{enumi}{0}
    \item \label{P2} $ \int_{\R^N} P(|y|)|x-y|^{2s-N}\,\dy \leq C|x|^{2s-N}$, for all $x \in \R^N \setminus B_1(0)$ and for some $C>0$.
\end{enumerate}
In what follows, we denote by  $(-\Delta)^{s}$ the fractional Laplace operator, defined as
$$
(-\Delta)^{s}u(x) := 2 \lim_{\epsilon \rightarrow 0^{+}} \int_{\mathbb{R}^{N} \backslash B_{\epsilon}(x)} \frac{u(x)-u(y)}{|x-y|^{N+2s}}\, \dy,
$$
where $s\in (0,1)$ and $B_{\epsilon}(x)$ is the ball of radius $\epsilon$ and centered at $x \in \mathbb{R}^{N}$. For $s=1$, problem \eqref{main-problem} reduces to the classical Laplacian case. In this setting, Yang et al.~\cite{Minbo-Abrantes-Pedro-Zhou-JDE-2025} studied the corresponding global problem and by combining the sub and supersolution method with variational techniques proved the existence of positive solutions. Moreover, they showed that for sufficiently large values of $\lambda$, the problem admits at least two distinct positive solutions. In the context of bounded domains, a seminal work was established by Brezis and Nirenberg~\cite{Brezis-1983}, who proved that local minimizers in $C^1$ are also local minimizers in the Sobolev space $H_0^1$. Building upon this result, Ambrosetti, Brezis, and Cerami~\cite{Brezis-Cerami-Ambrosetti-1994-JFA} obtained global results on existence, nonexistence and multiplicity of solutions depending on a real parameter. Later, Barrios et al.~\cite{Barrios-Colorado-Servadei-Soria-AIHP-2015} extended this framework to the fractional setting, proving the existence and multiplicity of positive solutions for critical problems involving the fractional Laplacian in bounded domains.

The objective of the present paper is to extend the results of Yang et al.~\cite{Minbo-Abrantes-Pedro-Zhou-JDE-2025} to the fractional setting, where the classical Laplacian is replaced by the fractional Laplacian $(-\Delta)^s$ with $s\in(0,1)$. Within this framework, we establish existence and non-existence results for positive solutions depending on the parameter $\lambda>0$. By means of an eigenvalue type problem (see Theorem \ref{Pezzo-Theorem-1.1} and Lemma \ref{Pezzo-Lemma-5.3}), combined with the sub and supersolution technique, we prove the existence of a critical value $\Lambda_s>0$ such that the qualitative behavior of positive solutions depends on the order relation between $\lambda$ and $\Lambda_s$. More precisely, we establish the following theorem.
\begin{theorem}[\textbf{Existence vs Nonexistence}]\label{Existence-first-solution}
For $s\in (0,1), N>2s$, assume that {\rm\ref{f1}-\ref{f3}, \ref{h}} and {\rm \ref{P1}} hold. Then there exists $\Lambda_s > 0$ as defined in \eqref{inf-A} such that,
\begin{enumerate}\rm
    \item[(i)] \eqref{main-problem} admits a positive solution $u_\lambda \in \mathcal{D}^{s,2}(\mathbb{R}^N)\cap L^\infty(\mathbb{R}^N)$ for all $\lambda > \Lambda_s$, and has no solution for all $\lambda < \Lambda_s$.
    \item[(ii)] Moreover, if $\lambda = \Lambda_s$ and \ref{f5} as given below holds, then \eqref{main-problem} admits at least one positive solution.
    \begin{enumerate}[label={($\bf f{5}$)}]
    \setcounter{enumi}{4}
    \item \label{f5}  $f$ is superlinear near $0$ i.e. $\displaystyle \lim\limits_{t \rightarrow 0^+}\frac{f(t)}{t}=0$.
   \end{enumerate}
\end{enumerate}
\end{theorem}
Under additional monotonicity and superlinearity assumptions on the function $f$ near the origin, we further prove that, for every $\lambda>\Lambda_s$, problem \eqref{main-problem} admits at least two distinct and ordered positive solutions. The first one is obtained by means of the sub and supersolution method, while the second is achieved through a classical linking result (see Theorem \ref{Linking-Theorem}), combined with our local minimization result of type $\mathcal{D}^{s,2}(\mathbb{R}^N)$ versus $C(\mathbb{R}^N;1+|x|^{N-2s})$ (see Theorem \ref{theorem-4.2}).

\begin{theorem}[\bf Multiplicity]\label{second-solution}
Let $N>2s,s\in (0,1)$. Assume that {\rm\ref{f1}-\ref{f3}, \ref{f5}, \ref{h}, \rm\ref{P1}} and {\rm\ref{P2}} hold. Then, for all $\lambda > \Lambda_s$, the problem \eqref{main-problem} admits at least two distinct positive solutions $v_\lambda, \widetilde{v}_\lambda \in \Ds2 \cap L^\infty (\R^N)$ satisfying $v_\lambda > \widetilde{v}_\lambda$ in $\R^N$.
\end{theorem}

\noindent\textbf{Organization of the paper.}
Section~\ref{Section-functional-framework} is devoted to establishing the functional framework and recalling some known results concerning a weighted eigenvalue problem. We also present certain compact embeddings of $\Ds2$, as well as regularity and positivity properties of the solutions. In Section~\ref{Section-Existence of solution}, we investigate the existence and nonexistence of solutions to \eqref{main-problem} by employing the sub- and supersolution method, depending on the range of the parameter~$\lambda$. Section~\ref{Section-Multiplicity} is concerned with the relation between the local minimizers of $\Ds2$ and those of the subspace~$X$ (see~\eqref{subspace-X} for the definition). Finally, we apply the Classical Linking Theorem to establish the existence of two distinct positive solutions to~\eqref{main-problem} when $\lambda$ is sufficiently large.
%%%%%%%%%%%%%%%%%%%%%%%%%%%%%%%%%%%%%%%%%%%%%%%%%%%%%%%%%%%%%%%%%%%%%%%%%%%%%%%%%%%%%%%%%%%%%%%%%%%%%%%%%%%%%
\section{Functional Framework and Some Technical Results}\label{Section-functional-framework}
In this section, we provide the functional framework associated with our problem \eqref{main-problem}. Moreover, we recall some known results of the corresponding weighted eigenvalue problem. Finally, we provide some embeddings of $\Ds2$ and prove the qualitative behavior of solutions, i.e., regularity and positivity.
\subsection{Functional Framework} Let $C_{c}^{\infty}(\R^N)$ be the vector space of all smooth functions with compact support. Recall that $\left(C_{c}^{\infty}(\R^N), \|\cdot\|_{\mathcal{D}}\right)$, where $\|\cdot\|_{\mathcal{D}}$ is given by \eqref{Gag-Norm}, is a normed vector space, but it is not complete. For $N>2s$ and $s\in (0,1)$, the completion of $C_{c}^{\infty}(\R^N)$ with respect to the norm $\|\cdot\|_{\mathcal{D}}$ is denoted by $\Ds2$ and is known as the fractional homogeneous Sobolev space. In other words, we write
\begin{align}\label{Gag-Norm}
    \Ds2 := \overline{C_{c}^{\infty}(\R^N)}^{\norm{\cdot}_{\mathcal{D}}}, \text{ where }  \|u\|_{\mathcal{D}} := \bigg(\Gag2\bigg)^{\frac{1}{2}}.
\end{align}
This fractional homogeneous Sobolev space has the following equivalent characterization (cf. \cite[Theorem 3.1]{Brasco-Castro-Vazquez-2021-CVPDE}):
\begin{align}\label{dsp}
   \Ds2 = \left\{ u \in L^{2_s^*}(\R^N) : \norm{u}_{\mathcal{D}} < \infty \right\}, \text{ where } 2_s^*= \frac{2N}{N-2s}. 
\end{align}
The above characterization also holds in the higher order case, that is, when $s \in (1,2)$ under the appropriate norm. The interested readers can refer to \cite{Biswas-Kumar-2024-higher} for the higher order characterization.

A function $u \in \Ds2$ is called a weak subsolution or supersolution respectively to \eqref{main-problem} if it satisfies
\begin{equation}\label{sub-super-weak-formulation}
        \iint_{\R^{2N}} \frac{( u(x)-u(y) )(\varphi(x)-\varphi(y))}{|x-y|^{N+2s}}\, \dxy 
   \leq (\text{ or }\geq ~) ~\lambda \int_{\R^N} h(x)f(u) \varphi(x) \, \dx, \text{ for all } \, \varphi \in \Ds2, \varphi \geq 0.
\end{equation}
If $u$ is both a weak subsolution and a supersolution, then it is called the weak solution to \eqref{main-problem}. In other words, $u \in \Ds2$ is said to be a weak solution to \eqref{main-problem} if it satisfies the following weak formulation i.e.
\begin{equation}\label{weak-formulation}
        \iint_{\R^{2N}} \frac{( u(x)-u(y) )(\varphi(x)-\varphi(y))}{|x-y|^{N+2s}}\, \dxy 
   =\lambda \int_{\R^N} h(x)f(u) \varphi(x) \, \dx, \text{ for all } \, \varphi \in \Ds2.
\end{equation}
The energy functional associated with \eqref{weak-formulation} is defined as $\mathcal{J}_\lambda: \Ds2 \rightarrow \mathbb{R}$ and given by
\begin{align}\label{energy-functional}
    \mathcal{J}_\lambda(u) = \frac{1}{2}\Gag2 - \lambda \int_{\R^N} h(x)F(u)\, \dx,\quad \text{ where } F(t) = \int_{0}^{t}f(s)\,\ds.
\end{align}
The energy functional $\mathcal{J}_\lambda$ is continuously Fr\'{e}chet differentiable on $\Ds2$ and its Fr\'{e}chet derivative is given by
\begin{align*}
    \langle \mathcal{J}_\lambda'(u), \varphi \rangle = \iint_{\mathbb{R}^N \times \mathbb{R}^N} \frac{(u(x)-u(y))(\varphi(x)-\varphi(y))}{|x - y|^{N + 2s}} \, \dxy - \lambda \int_{\mathbb{R}^N} h(x) f(u) \varphi(x) \, \dx, \text{ for all }\varphi \in \mathcal{D}^{s,2}(\mathbb{R}^N).
\end{align*}
Since $\mathcal{J}_\lambda'(u)$ belongs to the dual space of $\Ds2$, the above product $\langle \cdot, \cdot \rangle$ is indeed a duality product between the elements of the dual space of $\Ds2$ and the elements of $\Ds2$. Observe that every critical point of $\mathcal{J}_\lambda$ satisfies the weak formulation \eqref{weak-formulation}. Hence, every critical point of $\mathcal{J}_\lambda$ is also a weak solution to \eqref{main-problem}. Our aim is to show the existence of at least two critical points using variational arguments.

\subsection{Weighted Eigenvalue Problem} For $s \in (0,1)$ and $N>2s$, we consider the weighted eigenvalue problem
\begin{align}\label{weighted-problem}
    \begin{cases}
(-\Delta)^s u = \lambda h(x) u & \text{in } \mathbb{R}^N, \\
u \in \mathcal{D}^{s,2}(\mathbb{R}^N). \tag{$E_\lambda^s$}
\end{cases}
\end{align}
A pair $(u,\lambda) \in \mathcal{D}^{s,2}(\mathbb{R}^N) \times \mathbb{R}$ is a \textit{weak solution} of \eqref{weighted-problem} if the following Euler-Lagrange equation holds
\begin{align*}
    \iint_{\mathbb{R}^N \times \mathbb{R}^N} \frac{(u(x) - u(y))(\varphi(x) - \varphi(y))}{|x - y|^{N+2s}} \,\dxy = \lambda \int_{\mathbb{R}^N} h(x)u(x)\varphi(x)\,\dx, \text{ for all } \varphi \in \Ds2.
\end{align*}
The parameter $\lambda$ is called the eigenvalue and the corresponding function $u$ is called the eigenfunction of \eqref{weighted-problem} associated with $\lambda$.
\begin{theorem}{\rm(cf.\cite[Theorem 1.1]{Pezzo-Quaas-2020-NA})}\label{Pezzo-Theorem-1.1} Let $N>2s, s\in (0,1)$ and $h$ satisfy {\rm\ref{h}}. Then there exists a sequence $\left\{(u_n, \lambda_n(h)): n\in \mathbb{N} \right\}$ of eigenvalues and corresponding eigenfunctions of \eqref{weighted-problem} such that for every $n \in \mathbb{N}$,
\begin{align*}
\int_{\mathbb{R}^N} h(x)|u_n(x)|^2\, \dx = 1 \quad \text{ and } \quad 0 < \lambda_1(h) < \lambda_2(h) \leq \cdots \leq \lambda_n(h) \to \infty.
\end{align*}
The eigenvalue $\lambda_1(h)$ is called the first eigenvalue of \eqref{weighted-problem} and the corresponding eigenfunction is called the first eigenfunction. The eigenvalue $\lambda_1(h)$ is simple and the corresponding eigenfunctions do not change sign. Furthermore, $\lambda_1(h)$ admits the following characterization:
\begin{align*}
\lambda_1(h) = \min \left\{ \|u\|_{\mathcal{D}}^2 : u \in \mathcal{D}^{s,2}(\mathbb{R}^N) \text{  and  } \int_{\mathbb{R}^N} h(x)|u(x)|^2\, \dx = 1 \right\}.
\end{align*}    
\end{theorem}
\begin{lemma}{\rm(cf.\cite[Lemma 5.3]{Pezzo-Quaas-2020-NA})}\label{Pezzo-Lemma-5.3}
Let $N>2s, s\in (0,1)$ and $h$ satisfy {\rm\ref{h}}. Let $\lambda$ be an eigenvalue of \eqref{weighted-problem} and let $u$ be the corresponding eigenfunction. Then $u \in L^\infty(\mathbb{R}^N) \cap C^{0,\gamma}(\mathbb{R}^N)$ for some $\gamma \in (0,1)$.
\end{lemma}
Notice that Theorem \ref{Pezzo-Theorem-1.1} gives the existence of the first weighted eigenvalue $\lambda_1(h)$ of \eqref{weighted-problem} and it is characterize as follows:
\begin{align*}
\lambda_1(h) = \min_{u \in \mathcal{D}^{s,2}(\mathbb{R}^N) \setminus \{0\}} \frac{\displaystyle \iint_{\mathbb{R}^N \times \mathbb{R}^N} \frac{|u(x) - u(y)|^2}{|x - y|^{N + 2s}} \dxy}{\displaystyle \int_{\mathbb{R}^N} h(x)|u(x)|^2 \,\dx} > 0.
\end{align*}
Assume that $u_h \in \Ds2$ is the eigenfunction associated with $\lambda_1(h)$. Then, Lemma \ref{Pezzo-Lemma-5.3} infers that $\varphi_h$ is bounded and H\"{o}lder continuous with some exponent $\gamma \in (0,1)$.

\subsection{Some Embeddings, Regularity and Positivity of Solutions}
For a given nonnegative Lebesgue measurable function $g$, we define the space $L^2_g(\mathbb{R}^N)$ as the class of all real-valued Lebesgue measurable functions $u$ which satisfies
\begin{align*}
\int_{\mathbb{R}^N} g(x)|u(x)|^2 \dx < \infty.
\end{align*}
The space $L^2_g(\mathbb{R}^N)$ admits a real Hilbert space structure and its norm is induced by the following inner product.
\begin{align*}
\langle u, v \rangle_g := \int_{\mathbb{R}^N} g(x) u(x) v(x) \, \dx, \quad \text{for } u,v \in L^2_g(\mathbb{R}^N).
\end{align*}
It is important to mention that the embedding $\Ds2 \hookrightarrow L^{2_s^*}(\R^N)$ is continuous but it is not compact due to the lack of compactness in the fractional critical Sobolev embedding (cf.\cite[Theorem 6.5]{Hitchiker-guide-2012}). Our next result shows that the space $\Ds2$ is locally embedded into the Lebesgue spaces $L^q$ when $q$ is less than the critical exponent $2_s^*$. In addition, we establish the continuous and compact embedding of $\Ds2$ into some weighted Lebesgue spaces.
\begin{proposition}{\rm(cf. \cite[Proposition 2.2]{Biswas-Kumar-2024-Semipositone})}\label{embedding}
Let $N>2s$ and $q \in [1, 2_s^*)$. Then we have 
\begin{itemize}\rm
    \item[(i)] The embedding $\mathcal{D}^{s,2}(\mathbb{R}^N) \hookrightarrow L_{\text{loc}}^q(\mathbb{R}^N)$ is continuous and compact.
    \item[(ii)] For $0 \leq g \in L^\alpha(\mathbb{R}^N)$ with $\alpha = \frac{2_s^*}{2_s^* - q}$, the embedding $\mathcal{D}^{s,2}(\mathbb{R}^N) \hookrightarrow L_g^q(\mathbb{R}^N)$ is continuous and compact.
\end{itemize}
\end{proposition}
Since $h$ satisfies \ref{h}, it is easy to verify using the interpolation argument that $h \in L^\alpha(\mathbb{R}^N)$ for $\alpha = \frac{2_s^*}{2_s^* - q}$. Consequently, Proposition \ref{embedding} also holds under the assumption that $h$ satisfies \ref{h}.

In order to prove the regularity of solutions, we first define the tail space and local Sobolev spaces as follows:
\begin{align*}
L^{1}_{2s}(\mathbb{R}^N) = \left\{ u \in L^1_{\text{loc}}(\mathbb{R}^N) : \int_{\mathbb{R}^N} \frac{|u(x)|}{1 + |x|^{N + 2s}} \, \dx < \infty \right\},
\end{align*}
\begin{align*}
W^{s,2}_{\text{loc}}(\Omega) = \left\{ u \in L^2_{\text{loc}}(\Omega) :  \iint_{K \times K} \frac{|u(x) - u(y)|^2}{|x - y|^{N + 2s}} \, \dxy < \infty, \text{ for any compact set } K \subset \Omega \right\}.
\end{align*}
\begin{theorem}{\rm(cf.\cite[Theorem 1.1]{Nowak-2021-CVPDE})} \label{Nowak-Theorem-1.1}
Let $\Omega \subset \R^N$ be a domain, $N > 2s, s \in (0,1)$ and $g \in L^{q}_{\text{loc}}(\Omega)$ for some $q > \frac{N}{2s}$. Assume that $u \in W^{s,2}_{\text{loc}}(\Omega) \cap L^{1}_{2s}(\mathbb{R}^N)$ is a weak solution of $(-\Delta)^s u = g \text{ in } \Omega$ i.e., $u$ satisfies
\begin{align*}
    \iint_{\R^N \times \R^N} \frac{(u(x)-u(y))(\varphi(x)-\varphi(y))}{|x-y|^{N+2s}}\,\dxy = \int_{\Omega} g(x) \varphi(x)\,\dx, 
\end{align*}
for all $\varphi \in W^{s,2}(\Omega)$ such that $\varphi$ is compactly supported in $\Omega$. Then $u \in C^{0,\alpha}_{\text{loc}}(\Omega)$ for $\alpha \in \left(0, \min\left\{2s - \frac{N}{q}, 1\right\} \right)$. 
\end{theorem}
Further, we state the following technical lemma which will be useful in the subsequent proposition.
\begin{lemma}{\rm(cf.\cite[Lemma 3.1]{Iannizzotto-Mosconi-Squassina-2015-NoDEA})}\label{Iannizzotto-Lemma-3.1}
    For all $a, b \in \mathbb{R}$, $r \geq 2$ and $k > 0$, we have
\begin{align*}
(a - b)\left(a |a|_k^{r - 2} - b |b|_k^{r - 2} \right) \geq \frac{4(r - 1)}{r^2} \left( a |a|_k^{\frac{r}{2} - 1} - b |b|_k^{\frac{r}{2} - 1} \right)^2,
\end{align*}
where $t_k = \text{sgn}(t) \min \{ |t|, k \}$, $\text{sgn}(\cdot)$ is the signum function, and $|t|_k = \min \{ |t|, k \}$ for $k > 0, \, t \in \mathbb{R}$.
\end{lemma}
\begin{proposition}[\textbf{Regularity}]\label{regularity}
Let $N>2s, s\in (0,1)$, and {\rm\ref{f1}-- \ref{f3}} and {\rm\ref{h}} hold. Further, let $u \in \Ds2$ be a weak solution to \eqref{main-problem}. Then $u \in L^\infty(\mathbb{R}^N) \cap C(\mathbb{R}^N)$ and it satisfies the estimate
\begin{align*}
    \|u\|_{\infty} \leq C \|u\|_{2^*_s}, \text{ where } C=C(N,s,\lambda,C_0,\|h\|_\infty).
\end{align*}
\end{proposition}
\begin{proof}
Since $u \in \mathcal{D}^{s,2}(\mathbb{R}^N)$ is a weak solution to \eqref{main-problem}, it satisfies the following weak formulation
\begin{align}\label{reg-eq-1}
    \iint_{\mathbb{R}^N \times \mathbb{R}^N} \frac{(u(x) - u(y))(\varphi(x) - \varphi(y))}{|x - y|^{N + 2s}} \, \dxy = \lambda \int_{\mathbb{R}^N} h(x) f(u) \varphi(x) \, \dx, \text{ for all } \varphi \in \mathcal{D}^{s,2}(\mathbb{R}^N).
\end{align}
For $\tau>0$, we define a truncated function of $u$ as $u_\tau = \text{sgn}(u) \min \{ |u|, \tau \}$. For $r \geq 2$ and $k>0$, we set $\varphi_1 = u |u|_\tau^{r-2}$. Since the map $t \rightarrow t|t|_\tau^{r-2}$ is Lipschitz in $\R$ for all $r \geq 2$ and $k>0$, we get $\varphi_1 \in \mathcal{D}^{s,2}(\mathbb{R}^N)$. Taking $\varphi = \varphi_1$ into \eqref{reg-eq-1}, we obtain
\begin{align}\label{reg-eq-2}
    \iint_{\mathbb{R}^N \times \mathbb{R}^N} \frac{(u(x) - u(y)) \left(u(x) |u(x)|_\tau^{r-2} - u(y) |u(y)|_\tau^{r-2} \right)}{|x - y|^{N + 2s}} \, \dxy
= \lambda \int_{\mathbb{R}^N} h(x) f(u) u(x) |u(x)|_\tau^{r-2} \, \dx.
\end{align}
From Lemma \ref{Iannizzotto-Lemma-3.1} and the embedding $\mathcal{D}^{s,2}(\mathbb{R}^N) \hookrightarrow L^{2^*_s}(\mathbb{R}^N)$, we obtain the following estimates
\begin{align}\label{eq-A2}
    &\iint_{\mathbb{R}^N \times \mathbb{R}^N} \frac{(u(x) - u(y)) \left( u(x) |u(x)|_\tau^{r-2} - u(y) |u(y)|_\tau^{r-2} \right)}{|x - y|^{N + 2s}} \, \dxy \notag\\
    &\geq \frac{4(r - 1)}{r^2} \iint_{\mathbb{R}^N \times \mathbb{R}^N} \frac{\left( u(x) |u(x)|_\tau^{\frac{r}{2} - 1} - u(y) |u(y)|_\tau^{\frac{r}{2} - 1} \right)^2}{|x - y|^{N + 2s}} \, \dxy \notag\\
    &\geq \frac{4(r - 1)}{r^2} \, C(N,s) \left( \int_{\mathbb{R}^N} \left| u(x) |u(x)|_\tau^{\frac{r}{2} - 1} \right|^{2^*_s}\, \dx \right)^{2 / 2^*_s}.
\end{align}
Now we combine the identity \eqref{reg-eq-2} with the estimate \eqref{eq-A2} to get
\begin{align*}
\left( \int_{\mathbb{R}^N} \left| u(x) |u(x)|_\tau^{\frac{r}{2} - 1} \right|^{2^*_s}\, \dx \right)^{2 / 2^*_s}
\leq \frac{r^2}{4(r - 1)} \, C(N,s)^{-1} \, \lambda \int_{\mathbb{R}^N} h(x) f(u) |u(x)| |u(x)|_\tau^{r-2}\, \dx.
\end{align*}
Letting $\tau \to \infty$ and using the monotone convergence theorem yields
\begin{align}\label{eq-A3}
    \left( \int_{\mathbb{R}^N} |u(x)|^{\frac{r}{2}  2^*_s}\, \dx \right)^{\frac{2}{2^*_s}} 
\leq \frac{r^2}{4(r - 1)} \, C(N, s)^{-1} \, \lambda \int_{\mathbb{R}^N} h(x) f(u) |u(x)|^{r - 1}\, \dx.
\end{align}
\noindent \textbf{Step 1:} For $r = 2^*_s$, we show that $|u|^{r} \in L^{\frac{2^*_s}{2}}(\mathbb{R}^N)$. From \ref{f1} and \ref{f4}, we get $f(u) \leq C_0 |u|$, for every $u \in \Ds2$. Hence, we deduce that
\begin{align}\label{eq-A4}
    \int_{\mathbb{R}^N} h(x) f(u) |u(x)|^{r-1} \,\dx 
\leq C_0\int_{\mathbb{R}^N} h(x) |u(x)|^{r} \,\dx 
= C_0\int_{\mathbb{R}^N} h(x) |u(x)|^{2^*_s} \,\dx 
\leq C_0\|h\|_{\infty} \|u\|_{2^*_s}^{2^*_s}.
\end{align}
From \eqref{eq-A3} and \eqref{eq-A4}, we conclude that $|u|^{r} \in L^{\frac{2^*_s}{2}}(\mathbb{R}^N)$ for $r = 2^*_s$. 

\noindent
\textbf{Step 2.} In this step, we consider $r > 2^*_s$. Using \eqref{eq-A3} and \eqref{eq-A4}, we also write:
\begin{align}\label{eq-A5}
   \left( \int_{\mathbb{R}^N} |u(x)|^{\frac{r}{2}  2^*_s} \,\dx \right)^{\frac{2}{2^*_s}} \leq \frac{r^2}{4(r - 1)} \, C(N, s)^{-1} \lambda C_0 \, \|h\|_{\infty} \int_{\mathbb{R}^N} |u(x)|^r \,\dx.
\end{align}
We define a sequence of real numbers inductively as follows: $r_1 = 2^*_s$ and $r_{j+1} = \frac{2^*_s}{2} r_j$ for all $j \geq 1$. Then
\begin{align*}
   r_2 = \frac{2^*_s}{2} r_1, \quad 
r_3 = \left( \frac{2^*_s}{2} \right)^2 r_1, \quad 
\ldots, \quad 
r_{j+1} = \left( \frac{2^*_s}{2} \right)^j r_1 \quad \text{for all } j \geq 1. 
\end{align*}
Taking $r=r_j$ into \eqref{eq-A5}, we get
\begin{align*}
   \|u\|_{{r_{j+1}}}^{r_{j}} = \left( \int_{\mathbb{R}^N} |u(x)|^{r_{j+1}} \,\dx \right)^{\frac{2}{2^*_s}} 
\leq \frac{r_j^2}{4(r_j - 1)} \, C(N, s)^{-1} \lambda C_0 \|h\|_{\infty} \|u\|_{{r_j}}^{r_j},
\end{align*}
which further gives
\begin{align*}
    \|u\|_{{r_{j+1}}} \leq \left( \frac{r_j^2}{4(r_j - 1)}  C(N, s)^{-1} \lambda C_0 \|h\|_{\infty} \right)^{\frac{1}{r_{j}}} \|u\|_{{r_j}}.
\end{align*}
Since $r_j>2$ for all $j$, we get $\frac{r_j}{2(r_j - 1)} < 1$ and which further yields $\|u\|_{r_{j+1}} \leq (r_j C)^{1/r_{j}} \|u\|_{r_j}$, where $C=C(N, s)^{-1} \lambda C_0 \|h\|_{\infty}$. Upon iterating, we get
\begin{align}\label{eq-A7}
    \|u\|_{r_{j+1}} \leq (r_{j} C)^{1/r_{j}} \cdots (r_1 C)^{1/r_1} \|u\|_{r_1}
= \left( \prod_{k=1}^{j} (r_k C)^{1/r_k} \right) \|u\|_{2_s^*}.
\end{align}
Notice that
\begin{align*}
    \sum_{k=1}^{\infty} \frac{1}{r_k} = \frac{1}{r_1}+\sum_{k=1}^{\infty} \frac{1}{r_{k+1}} = \frac{1}{2_s^*} +\frac{1}{2_s^*} \sum_{k=1}^{\infty} \left(\frac{2}{2_s^*} \right)^k = \frac{1}{2_s^*-2},
\end{align*} and
\begin{align}\label{product-rk}
    \prod_{k=1}^{\infty} (r_k)^{1/r_k} = \exp\left(\sum_{k=1}^{\infty} \frac{\log(r_k)}{r_k} \right)= \exp\left( \frac{\log(r_1)}{r_1} +\sum_{k=1}^{\infty} \frac{\log(r_{k+1})}{r_{k+1}} \right).
\end{align}
We can write the second series given above as
\begin{align*}
    \sum_{k=1}^{\infty} \frac{\log(r_{k+1})}{r_{k+1}} = \sum_{k=1}^{\infty} \frac{k \log(\frac{2_s^*}{2})+ \log(2_s^*)}{\left( \frac{2_s^*}{2}\right)^k 2_s^*} = \frac{1}{2_s^*}\log\left(\frac{2_s^*}{2}\right)\sum_{k=1}^{\infty} \left( \frac{2}{2_s^*}\right)^k k + \frac{\log(2_s^*)}{2_s^*}\sum_{k=1}^{\infty} \left( \frac{2}{2_s^*}\right)^k.
\end{align*}
The two series on the right hand side are convergent. Thus, the product \eqref{product-rk} is finite. Now we take the limit as $j \rightarrow \infty$ in \eqref{eq-A7} and get $\|u\|_{\infty} \leq C \|u\|_{2^*_s}$, where $C$ is independent of $r_j$. Consequently, we get $u \in L^\infty(\mathbb{R}^N)$.

Next we claim that $u \in C(\R^N)$. Assume that $\Omega \subset \R^N$ is an open bounded Lipschitz set and $ \varphi_2 \in W^{s,2}(\Omega)$ such that $\varphi_2$ is compactly supported in $\Omega$. Define $\tilde{\varphi}_2(x) = \varphi_2(x)$ for $x \in \Omega$, and $\tilde{\varphi}_2(x) = 0$ for $x \in \R^N \setminus \Omega$. Then, $\tilde{\varphi}_2 \in W^{s,2}(\R^N)$ due to \cite[Theorem 5.1]{Hitchiker-guide-2012}. The fractional Sobolev embedding infers that $\tilde{\varphi}_2 \in L^{2_s^*} (\R^N)$. Thus, $\tilde{\varphi}_2 \in \Ds2$ (see the definition of $\Ds2$ given in \eqref{dsp}). Since $u \in \Ds2 \cap L^{\infty}(\R^N)$, the embeddings $\Ds2 \hookrightarrow L^{2^*_s}(\R^N)$ and $L^{2^*_s}(\Omega) \hookrightarrow L^2(\Omega)$ infer $\Ds2 \subset W^{s,2}_{\text{loc}}(\Omega)$. Note that $L^{\infty}(\R^N) \subset L^{1}_{2s}(\R^N)$. Therefore, in view of \eqref{reg-eq-1}, $u \in W^{s,2}_{\text{loc}}(\Omega) \cap L^{1}_{2s}(\R^N)$ satisfies the following identity:
\begin{align*}
    \iint_{\R^{2N}} \frac{\big(u(x)-u(y)\big)\big(\tilde{\varphi}_2(x) - \tilde{\varphi}_2(y)\big)}{|x-y|^{N+sp}}\, \dxy & = \lambda \int_{\R^N }  h(x) f(u) \tilde{\varphi}_2(x) \, \dx, \\
    & =\lambda \int_{ \Omega } h(x) f(u) \varphi_2(x) \, \dx, ~~ \forall \varphi_2 \in W^{s,2}(\Omega), \text{supp}(\varphi_2) \subset \Omega.
\end{align*}
For $q > \frac{N}{2s}$, we notice that
\begin{align*}
    \int_{\R^N} \left(h(x) f(u)\right)^q \,\dx \leq C_0^q \int_{\R^N} h(x)^q |u|^q\,\dx \leq C_0^q \|u\|_\infty^q \|h\|_q^q.
\end{align*}
The first inequality follows using $\abs{f(u)} \leq C_0 \abs{u}$, for all $u \in \Ds2$. The quantity on the right hand side of the above inequality is finite using \ref{h}, which implies $hf(u) \in L^q(\R^N)$ for $q > \frac{N}{2s}$. Thus, $u$ solves weakly the equation
\[ (-\Delta)_{p}^{s} u = \lambda h(x) f(u) \text{ in } \Omega.\]
We apply Theorem \ref{Nowak-Theorem-1.1} to get $u \in C_{\text{loc}}^{0,\alpha}(\Omega)$ for some $\alpha \in (0,1)$. In particular, $u \in C_{\text{loc}}(\Omega)$ which yields $u \in C(\Omega)$ for every bounded open Lipschitz set $\Omega \subset \R^N$. Next, for any compact set $K \subset \R^N$, we have $K \subset \Omega$ for some bounded open Lipschitz set $\Omega$. Thus, $u \in C(K)$ as well which implies $u \in C_{\text{loc}}(\R^N)$. Hence,  $u \in C(\R^N)$.
\end{proof}
In the next result, we show that every nontrivial weak solution to \eqref{main-problem} is indeed positive.
\begin{proposition}[\textbf{Positivity}]\label{Positivity}
   Let {\rm\ref{f1}--\rm\ref{f3}} and {\rm\ref{h}} hold. Then the nontrivial critical points of $\mathcal{J}_\lambda$ are positive.
\end{proposition}
\begin{proof}
Our first step is to show that any nontrivial solution to \eqref{main-problem} is nonnegative. In the second step, we will apply the strong maximum principle to obtain the positivity. Let $u$ be a nontrivial solution to \eqref{main-problem}. Let $A=\{x \in \R^N : u(x) \geq 0\}$ be the set where $u$ becomes nonnegative. Then, for every $\varphi \in \Ds2$ we have
\begin{align}\label{Id1}
        \iint_{\R^{2N}} \frac{( u(x)-u(y) )(\varphi(x)-\varphi(y))}{|x-y|^{N+2s}}\, \dxy 
   =\lambda \int_{\R^N} h(x)f(u) \varphi(x) \, \dx=\lambda \int_{A} h(x)f(u) \varphi(x) \, \dx.
\end{align}
The last identity follows since $f(t)= 0$ for all $t \leq 0$. Let $u^+=\max\{u,0\}$ and $u^- = \max\{-u,0\}$. Taking $\varphi = - u^-$ in the above identity, the right hand side integral becomes zero since $u^- \equiv 0$ on $A$. Moreover, we have
\begin{align*}
    (u(x)-u(y))(\varphi(x)-\varphi(y)) &= -(u(x)-u(y))(u^-(x)- u^-(y))\\
    &= -(u^+(x)-u^+(y))(u^-(x)- u^-(y)) + (u^-(x)- u^-(y))^2\\
    &=u^+(x)u^-(y) + u^+(y) u^-(x) + (u^-(x)- u^-(y))^2 \\
    & \geq (u^-(x)- u^-(y))^2.
\end{align*}
From \eqref{Id1} and the above estimate, we get
\begin{align*}
  \iint_{\R^{2N}} \frac{(u^-(x)- u^-(y))^2}{|x-y|^{N+2s}}\, \dxy \leq   \iint_{\R^{2N}} \frac{( u(x)-u(y) )(\varphi(x)-\varphi(y))}{|x-y|^{N+2s}}\, \dxy =0,
\end{align*}
which further implies $\|u^-\|_{\mathcal{D}}=0$. The Sobolev embedding $\Ds2 \hookrightarrow L^{2_s^*}(\R^N)$ infers $u^- = 0$ a.e. in $\R^N$. Subsequently, $u \geq 0$ a.e. in $\R^N$. Indeed, $u \geq 0$ in $\R^N$ since $u$ is continuous (cf. Proposition \ref{regularity}). To prove positivity, we assume on contrary that there exists some $ x_0 \in \mathbb{R}^N$ such that $u(x_0) = 0$. The strong maximum principle \cite[Proposition 5.2.1]{Dipierro-Medina-Valdinoci-2017-Book} gives $u \equiv 0$, which is a contradiction. Hence, the result follows immediately.
\end{proof}

\section{Existence of a weak solution}\label{Section-Existence of solution}
In this section, we define a critical threshold level $\Lambda_s$ for the parameter $\lambda$ and prove the existence and nonexistence results depending on the order relation between $\lambda$ and $\Lambda_s$.
\begin{proposition}\label{Existence-1}
Let {\rm \ref{f1}, \ref{f2}} and {\rm\ref{h}} hold. Then, there exists $\lambda_0 > 0$ such that \eqref{main-problem} has a positive solution whenever $\lambda \geq \lambda_0$.
\end{proposition}
\begin{proof}
For any $\varphi \in C_c^\infty(\mathbb{R}^N)$ with $\varphi > 0$, we have $\mathcal{J}_\lambda(\varphi)<0$ when we take $\lambda$ sufficiently large. Fixing $\lambda_0 > 0$ such that $\mathcal{J}_\lambda(\varphi)<0$ whenever $\lambda \geq \lambda_0$. Using \ref{f2}, there exists a real number $t_\lambda > 0$ such that
\begin{align*}
f(t) \leq \frac{\lambda_1(h)}{2\lambda} t, \quad \text{ for all } t \geq t_\lambda,
\end{align*}
where $\lambda_1(h)$ is the first eigenvalue of \eqref{weighted-problem}. Thus, $F(t) \leq \frac{\lambda_1(h)}{4\lambda} t^2 + C_\lambda$, for all $t \in \R$, where $C_\lambda$ is a constant. Now
\begin{align}\label{Exis-1-eq-1}
    \mathcal{J}_\lambda(u) &= \frac{1}{2} \iint_{\mathbb{R}^N \times \mathbb{R}^N} \frac{|u(x)-u(y)|^2}{|x-y|^{N+2s}}\, \dxy - \lambda \int_{\mathbb{R}^N} h(x) F(u)\, \dx \notag\\
    & \geq \frac{1}{2} \iint_{\mathbb{R}^N \times \mathbb{R}^N} \frac{|u(x)-u(y)|^2}{|x-y|^{N+2s}}\, \dxy - \frac{\lambda_1(h)}{4} \int_{\mathbb{R}^N} h(x) u^2\, \dx - \lambda C_\lambda\int_{\mathbb{R}^N} h(x)\,\dx \notag\\
    &\geq \frac{1}{4} \iint_{\mathbb{R}^N \times \mathbb{R}^N} \frac{|u(x)-u(y)|^2}{|x-y|^{N+2s}}\, \dxy - \lambda C_\lambda \|h\|_1.
\end{align}
From \eqref{Exis-1-eq-1}, we obtain that $\mathcal{J}_\lambda$ is bounded from below on $\Ds2$. Let $\|\varphi\|_{\mathcal{D}}< r_\lambda$, for some $r_\lambda>0$. Then $\varphi \in B_{r_\lambda}(0) \subset \Ds2$ and $m_\lambda = \inf \left\{\mathcal{J}_\lambda(v): v \in \overline{B}_{r_\lambda}(0)\right\} \leq \mathcal{J}_\lambda(\varphi) <0$, for all $\lambda \geq \lambda_0$. Since $\mathcal{J}_\lambda$ is a weakly lower semi-continuous function on $\Ds2$, it attains its minimum on the ball $\overline{B}_{r_\lambda}(0)$ i.e., there exists a function $u_\lambda \in \overline{B}_{r_\lambda}(0)$ such that $m_\lambda = \mathcal{J}_\lambda(u_\lambda)$. If $u_\lambda \in \partial B_{r_\lambda}(0)$, then $\|u_\lambda\|=r_\lambda$ and we can choose $r_\lambda>0$ large enough such that $m_\lambda=\mathcal{J}_\lambda(u_\lambda)>0$ due to \eqref{Exis-1-eq-1}, which gives a contradiction to the fact that $m_\lambda<0$. Consequently, we get $u_\lambda \notin \partial B_{r_\lambda}(0)$ which implies $u_\lambda$ is a critical point of $\mathcal{J}_\lambda$ and it is positive by Proposition \ref{Positivity}. Hence, $u_\lambda$ is a positive solution to \eqref{main-problem} for $\lambda \geq \lambda_0,$ where $\lambda_0$ is sufficiently large.
\end{proof}
Let us define a set $\mathcal{A}_s$ as follows:
\begin{align}\label{set-A}
  \boldsymbol{  \mathcal{A}_s = \left\{ \lambda \in \mathbb{R}: \eqref{main-problem} \textbf{ admits a positive solution}\right\}}.
\end{align}
The set $\mathcal{A}_s$ is non-empty and bounded from below due to Proposition \ref{Existence-1}. Therefore, we define
\begin{align}\label{inf-A}
  \boldsymbol{ \Lambda_s = \inf\mathcal{A}_s}.
\end{align}
\begin{lemma}
Let {\rm\ref{f1}--\ref{f3}} and {\rm\ref{h}} hold. Then $\Lambda_s$ is positive.
\end{lemma}
\begin{proof}
For each $\lambda \in \mathcal{A}_s$, by \ref{f4} we have that the solution of \eqref{main-problem} verifies
\begin{align*}
\iint_{\mathbb{R}^N \times \mathbb{R}^N} \frac{|u(x)-u(y)|^2}{|x-y|^{N+2s}}\, \dxy = \lambda \int_{\mathbb{R}^N} h(x)f(u)u \, \dx \leq \lambda C_0 \int_{\mathbb{R}^N} h(x)u^2 \, \dx.
\end{align*}
Therefore
\begin{align*}
\lambda_1(h) \leq \frac{\iint_{\mathbb{R}^N \times \mathbb{R}^N} \frac{|u(x)-u(y)|^2}{|x-y|^{N+2s}}\, \dxy}{\int_{\mathbb{R}^N} h(x)u^2 \, \dx} \leq \lambda C_0,
\end{align*}
which implies
\begin{align*}
\lambda \geq \frac{\lambda_1(h)}{C_0} > 0.
\end{align*}
\end{proof}
In the following, we prove the existence and nonexistence of solutions to \eqref{main-problem} depending on the order relation between the parameters $\lambda$ and $\Lambda_s$. Following the proof of \cite[Proposition 3.4]{Fractional-Sublinear-2016-AA} and \cite[Lemma 3.5]{Fractional-Sublinear-2016-AA}, we deduce the following Theorem
\begin{theorem}\label{Theorem-3.1}
    If $N>2s$ and $h$ satisfies {\rm\ref{h}}, then there exists a unique bounded positive solution $u\in \Ds2$ satisfying
    \begin{align}\label{Linear-Problem}
\begin{cases}
        (-\Delta)^s u = h \text{ in } \R^N,\\
        \lim\limits_{|x| \rightarrow \infty} u(x)=0.
    \end{cases}
    \end{align}
\end{theorem}
% \tr{\begin{theorem}[Proposition 3.4, Lemma 3.5, \cite{Fractional-Sublinear-2016-AA}]
%     If $N>2s$ and $0 \leq h(x) \leq C (1+|x|^{-\beta})$ for all $x \in \R^N$, for some $C>0$ and $\beta \in (N,\infty)$, then $K_s * h \in \Ds2 \cap L^\infty(\R^N)$ and it is a weak solution to the  equation
%     \begin{align}
%     \begin{cases}
%         (-\Delta)^s u = h \text{ in } \R^N,\\
%         \lim\limits_{|x| \rightarrow \infty} u(x)=0.
%     \end{cases}
%     \end{align}
% The convolution $K_s * h$ is given by $(K_s * h)(x) = \int_{\R^N} \frac{h(y)}{|x-y|^{N-2s}}\,\dy$, for all $x\in \R^N$. Moreover, if $v \in \Ds2$ is any bounded solution to \eqref{Linear-Problem}, then it coincides with $K_s * h$, i.e., the solution is unique.
% \end{theorem}
% \begin{remark}\label{Remark-3.2}
%  Theorem \ref{Theorem-3.1} also holds under the assumptions \ref{h}, \ref{P1} and $N>4s$ (cf. \cite[Remark 3.3]{Fractional-Sublinear-2016-AA}).
% \end{remark}}

\begin{proof}[Proof of Theorem \ref{Existence-first-solution}]
(i) By the definition of $\Lambda_s$, for every $\epsilon>0$ there exists $\beta \in \mathcal{A}_s$ such that $\beta < \Lambda_s + \epsilon$. For $\lambda>\Lambda_s$, we choose $\beta \in (\Lambda_s,\lambda)$ such that $\beta \in \mathcal{A}_s$, i.e., there exists a positive solution $u_\beta \in \Ds2 \cap L^\infty (\R^N)$ to \eqref{main-problem} when $\lambda=\beta$. Moreover, $u_\beta$ is a positive subsolution to \eqref{main-problem}. The next step is to construct a positive supersolution to \eqref{main-problem}. From Theorem \ref{Theorem-3.1}, let $u_h\in \Ds2 \cap L^\infty(\R^N)$ be the unique positive bounded solution to \eqref{Linear-Problem}. From \ref{f2}, there exists $s_\lambda>0$ such that 
\begin{align*}
    \frac{f(s)}{s} \leq \frac{1}{\lambda \|u_h\|_\infty}, \text{ for all } s>s_\lambda.
\end{align*}
Since $f$ is a continuous function on $\R$, we define a new function $\widetilde{f}$ as follows:
\begin{align*}
    \widetilde{f}(t):=\max_{0\leq s \leq t} f(s), t \in [0,\infty).
\end{align*}
Observe that $\widetilde{f}(t) \geq f(t)$, for all $t \in [0,+\infty)$, and $\widetilde{f}$ is non-decreasing in $[0,+\infty)$. Let $r_\lambda > \max\left\{ \frac{s_\lambda}{\|u_h\|_\infty}, \lambda C_0 s_\lambda \right\}$, where $C_0$ is given in \ref{f4}. Then we have
\begin{align*}
    \frac{\widetilde{f}(r_\lambda \|u_h\|_\infty)}{r_\lambda \|u_h\|_\infty}& = \frac{\max\{f(s): 0\leq s \leq r_\lambda \|u_h\|_\infty \}}{r_\lambda \|u_h\|_\infty}\\
    &= \frac{\max\{f(s) \chi_{[0,s_\lambda]}(s) + f(s) \chi_{[s_\lambda, r_\lambda \|u_h\|_\infty]}(s) : 0\leq s \leq r_\lambda \|u_h\|_\infty \}}{r_\lambda \|u_h\|_\infty}\\
    &= \max\left\{ \frac{f(s) \chi_{[0,s_\lambda]}(s)}{r_\lambda \|u_h\|_\infty} + \frac{f(s) \chi_{[s_\lambda, r_\lambda \|u_h\|_\infty]}(s)}{r_\lambda \|u_h\|_\infty} : 0\leq s \leq r_\lambda \|u_h\|_\infty \right\}\\
    &\leq \max\left\{ \frac{C_0 s \chi_{[0,s_\lambda]}(s)}{r_\lambda \|u_h\|_\infty} + \frac{f(s)}{s} \frac{s \chi_{[s_\lambda, r_\lambda \|u_h\|_\infty]}(s)}{r_\lambda \|u_h\|_\infty} : 0\leq s \leq r_\lambda \|u_h\|_\infty \right\}\\
    &\leq \max\left\{ \frac{ \chi_{[0,s_\lambda]}(s)}{\lambda \|u_h\|_\infty} + \frac{ \chi_{[s_\lambda, r_\lambda \|u_h\|_\infty]}(s)}{\lambda \|u_h\|_\infty} : 0\leq s \leq r_\lambda \|u_h\|_\infty \right\}\\
    &=\frac{1}{\lambda \|u_h\|_\infty},
\end{align*}
where $\chi_A$ is the characteristic function on the set $A$. Let us define $v_\lambda (x) = r_\lambda u_h(x)$. Then $v_\lambda$ satisfies the following inequality on $\R^N$ in the weak sense
\begin{align}\label{Super-soluton}
    (-\Delta)^s v_\lambda &= r_\lambda (-\Delta)^s u_h = r_\lambda h(x)\notag\\
    & \geq \lambda \widetilde{f}(r_\lambda \|u_h\|_\infty) h(x) \notag \\
    & \geq \lambda h(x) \widetilde{f}(r_\lambda u_h) = \lambda h(x) \widetilde{f}(v_\lambda) \notag \\
    & \geq \lambda h(x) f(v_\lambda).
\end{align}
From \eqref{Super-soluton}, $v_\lambda$ is a positive super solution to \eqref{main-problem}. Recall that $u_\beta$ and $v_\lambda$ satisfy the weak formulations
\begin{equation}\label{u-beta-Euler}
        \iint_{\R^{2N}} \frac{( u_\beta(x)-u_\beta(y) )(\varphi(x)-\varphi(y))}{|x-y|^{N+2s}}\, \dxy 
   =\beta \int_{\R^N} h(x)f(u_\beta) \varphi(x) \, \dx, \quad \text{ for all } \varphi \in \Ds2,
\end{equation}
and 
\begin{equation}\label{v-lambda-Euler}
        \iint_{\R^{2N}} \frac{( v_\lambda(x)-v_\lambda(y) )(\varphi(x)-\varphi(y))}{|x-y|^{N+2s}}\, \dxy 
   =r_\lambda \int_{\R^N} h(x) \varphi(x) \, \dx, \quad \text{ for all } \varphi \in \Ds2.
\end{equation}
Now we fix $r_\lambda > \max\left\{\beta C_0 \|u_\beta \|_\infty, \frac{s_\lambda}{\|u_h\|_\infty}, \lambda C_0 s_\lambda \right\}$. Subtracting \eqref{v-lambda-Euler} from \eqref{u-beta-Euler} and taking $\varphi = (u_\beta - v_\lambda)^+$ as a test function, we get
\begin{align*}
    0 &\leq \iint_{\R^{2N}} \frac{| (u_\beta - v_\lambda)^+(x)-(u_\beta - v_\lambda)^+(y) |^2}{|x-y|^{N+2s}}\, \dxy \\
    & \leq  \iint_{\R^{2N}} \frac{ ((u_\beta - v_\lambda)(x)-(u_\beta - v_\lambda)(y)) ((u_\beta - v_\lambda)^+(x)-(u_\beta - v_\lambda)^+(y)) }{|x-y|^{N+2s}}\, \dxy\\
    &= \int_{\R^N} \left[ \beta f(u_\beta)- r_\lambda \right]  h(x) (u_\beta - v_\lambda)^+(x) \, \dx \\
    &\leq \int_{\R^N} \left[ \beta C_0 \|u_\beta\|_\infty - r_\lambda \right]  h(x) (u_\beta - v_\lambda)^+(x) \, \dx \\
    & \leq 0.
\end{align*}
The above chain of inequalities yields $\|(u_\beta - v_\lambda)^+\|_{\mathcal{D}}=0$, which further gives $(u_\beta - v_\lambda)^+ =0$ in $\R^N$. Therefore
\begin{align}\label{solution-compare}
    u_\beta (x) \leq v_\lambda (x), \text{ for all } x \in \R^N.
\end{align}
From the monotone iteration technique used in \cite[Lemma 5.4]{Fractional-Sublinear-2016-AA} and using the comparison principle given in \cite[Lemma 3.6]{Fractional-Sublinear-2016-AA}, we obtain the existence of a positive and bounded solution $u_\lambda$ to \eqref{main-problem} for all $\lambda > \Lambda_s$. Moreover, by the definition of $\Lambda_s$ as given in \eqref{inf-A}, it follows that \eqref{main-problem} has no solution for all $\lambda < \Lambda_s$.

(ii) Assume that $\lambda_n > \Lambda_s$ is a decreasing sequence converging to $\Lambda_s$ and $u_n$ is a solution to \eqref{main-problem} with $\lambda=\lambda_n$. Then we have
\begin{align*}
    \iint_{\R^N \times \R^N} \frac{(u_n(x)-u_n(y))(\varphi(x)-\varphi(y))}{|x-y|^{N+2s}}\,\dxy = \lambda_n \int_{\R^N} h(x)f(u_n)\varphi(x)\,\dx, \text{ for all }\varphi \in \Ds2.
\end{align*}
Using monotone iteration, we obtain a sequence of solutions given by $v_{\lambda_1} \geq u_1 \geq u_2 \cdots$, which implies
\begin{align*}
    \iint_{\R^N \times \R^N} \frac{|u_n(x)-u_n(y)|^2}{|x-y|^{N+2s}}\,\dxy \leq \lambda_1 C_0 \int_{\R^N} h(x)(u_1)^2\,\dx < +\infty.
\end{align*}
Clearly, $\{u_n\}$ is bounded in $\Ds2$. Consequently, $\exists$ $u_{\Lambda_s} \in \Ds2$ such that $0 \leq u_{\Lambda_s} \leq v_{\lambda_1}$ in $\R^N$ and
\begin{align*}
u_n \rightharpoonup u_{\Lambda_s} &\text{ in } \Ds2,\\
u_n \rightarrow u_{\Lambda_s} &\text{ a.e. in } \R^N.
\end{align*}
Observe that
\begin{align*}
    \iint_{\R^N \times \R^N} \frac{(u_n(x)-u_n(y))(\varphi(x)-\varphi(y))}{|x-y|^{N+2s}}\,\dxy \rightarrow \iint_{\R^N \times \R^N} \frac{(u_{\Lambda_s}(x)-u_{\Lambda_s}(y))(\varphi(x)-\varphi(y))}{|x-y|^{N+2s}}\,\dxy \text{ as } n \rightarrow \infty,
\end{align*}
and
\begin{align*}
    |h(x) [f(u_n)-f(u_{\Lambda_s})]\varphi| \leq C_1 h(x) |\varphi|, \text{ for all } n \in \mathbb{N}.
\end{align*}
Since $h$ satisfies \ref{h} and $\varphi \in \Ds2 \hookrightarrow L^{2_s^*}(\R^N)$, we get $h \varphi \in L^1(\R^N)$. Using dominated convergence theorem, we obtain
\begin{align*}
    \iint_{\R^N \times \R^N} \frac{(u_{\Lambda_s}(x)-u_{\Lambda_s}(y))(\varphi(x)-\varphi(y))}{|x-y|^{N+2s}}\,\dxy = {\Lambda_s} \int_{\R^N} h(x) f(u_{\Lambda_s}) \varphi (x)\,\dx, \text{ for all } \varphi \in \Ds2.
\end{align*}
Thus, $u_{\Lambda_s}$ is a weak solution to \eqref{main-problem} when $\lambda = {\Lambda_s}$. Next we claim that $u_{\Lambda_s}$ is positive. First we show that
\begin{align*}
    u_n \rightarrow u_{\Lambda_s} \text{ in } L^\infty(\R^N).
\end{align*}
Since $|x|^{2s-N}$ is a fundamental solution of $(-\Delta)^s$ (see \cite[Theorem 5]{Stinga-2019-Book}), by Riesz representation formula we get
\begin{align*}
u_n(x) = \lambda_n C(N,s) \int_{\mathbb{R}^N} \frac{h(y) f(u_n(y))}{|x-y|^{N-2s}} \, dy, \text{ for all } x \in \R^N \text{ and for all } n \in \mathbb{N}, 
\end{align*} and 
\begin{align*}
u_{\Lambda_s}(x) = {\Lambda_s} C(N,s) \int_{\mathbb{R}^N} \frac{h(y) f(u_{\Lambda_s}(y))}{|x-y|^{N-2s}} \, dy, \text{ for all } x \in \R^N,
\end{align*}
for some positive constant $C(N,s)$. Now we estimate $|u_n - u_{\Lambda_s}|$ as follows:
\begin{align}\label{sec3-eq1}
    |u_n(x) - u_{\Lambda_s}(x)| &\leq 
C(N,s)  \int_{B_1(x)} \frac{h(y) \big| \lambda_n f(u_n(y)) - {\Lambda_s} f(u_{\Lambda_s}(y)) \big|}{|x-y|^{N-2s}} \, \dy \notag\\
& \qquad+ C(N,s) \int_{\mathbb{R}^N \setminus B_1(x)} \frac{h(y) \big|\lambda_n f(u_n(y)) - {\Lambda_s} f(u_{\Lambda_s}(y)) \big|}{|x-y|^{N-2s}} \, \dy.
\end{align}
Take $\delta \in \left(1 , \frac{N}{N-2s}\right)$. Using H\"older’s inequality with conjugate pair $(\delta, \delta')$, we estimate the first integral of \eqref{sec3-eq1}, i.e.,
\begin{align}\label{sec3-eq2}
    &\int_{B_1(x)} \frac{h(y) \big|\lambda_n f(u_n(y)) -{\Lambda_s} f(u_{\Lambda_s}(y)) \big|}{|x-y|^{N-2s}} \, \dy \notag\\
    & \leq \left( \int_{B_1(x)} \frac{1}{|x-y|^{(N-2s)\delta}} \, \dy \right)^{\tfrac{1}{\delta}} 
\left( \int_{B_1(x)} h(y)^{\delta'} \big| \lambda_n f(u_n(y)) - {\Lambda_s} f(u_{\Lambda_s}(y)) \big|^{\delta'} \, \dy \right)^{\tfrac{1}{\delta'}}.
\end{align}
We calculate 
\begin{align*}
\int_{B_1(x)} \frac{1}{|x-y|^{(N-2s)\delta}} \, dy 
= \omega_N \int_0^1 r^{N-1} r^{-(N-2s)\delta} \, dr \leq C_1(N,s).
\end{align*}
Now we prove that the second integral in \eqref{sec3-eq2} converges to zero. Observe that
\begin{align*}
    \big|\lambda_n f(u_n(y)) - {\Lambda_s} f(u_{\Lambda_s}(y)) \big|^{\delta'}& \leq 2^{\delta'-1} 
\left( (\lambda_n - {\Lambda_s})^{\delta'} f(u_n(y) )^{\delta'} +{\Lambda_s}^{\delta'} \big| f(u_n(y)) - f(u_{\Lambda_s}(y)) \big|^{\delta'} \right)\\
&\leq 2^{\delta'-1} 
\left( (\lambda_n - {\Lambda_s})^{\delta'} f(u_n(y) )^{\delta'} +{\Lambda_s}^{\delta'} L |u_n(y) - u_{\Lambda_s}(y)|^{\delta'} \right)\\
& \leq 2^{\delta'-1} 
\left( (\lambda_1 - {\Lambda_s})^{\delta'} (C_0 (\delta
') u_1(y) )^{\delta'} +{\Lambda_s}^{\delta'} L C(\delta') \left( |u_1(y)|^{\delta'} + |u_{\Lambda_s}(y)|^{\delta'} \right) \right)\\
& \leq C_3 (L,\delta', \lambda, {\Lambda_s}, \|u_1\|_\infty, \|u_{\Lambda_s}\|_\infty),
\end{align*}
where $L$ is the Lipschitz constant. Since, $f(u_n(y)) \to f(u_{\Lambda_s}(y))$ and $\lambda_n \to {\Lambda_s}$ as $n \to \infty$, we apply the dominated convergence theorem to get 
\begin{align}\label{sec3-eq3}
    \int_{B_1(x)} h(y)^{\delta'}  |\lambda_n f(u_n(y)) - {\Lambda_s} f(u_{\Lambda_s}(y))|^{\delta'}  \, \dy \to 0.
\end{align}
Next, the second integral of \eqref{sec3-eq1} has the following bound:
\begin{align*}
    \int_{\mathbb{R}^N \setminus B_1(x)} 
\frac{h(y) \left| \lambda_n f(u_n(y)) - {\Lambda_s} f(u_{\Lambda_s}(y)) \right|}{|x-y|^{N-2s}} \, \dy 
&\leq \int_{\mathbb{R}^N \setminus B_1(x)} h(y) \left| \lambda_n f(u_n(y)) - {\Lambda_s} f(u_{\Lambda_s}(y)) \right|\,\dy\\
& \leq \int_{\mathbb{R}^N} h(y) \left| \lambda_n f(u_n(y)) - {\Lambda_s} f(u_{\Lambda_s}(y)) \right|\,\dy.
\end{align*}
Again we apply the dominated convergence theorem and obtain
\begin{align}\label{sec3-eq4}
    \int_{\mathbb{R}^N} h(y) \left|\lambda_n f(u_n(y)) - {\Lambda_s} f(u_{\Lambda_s}(y)) \right| \to 0.
\end{align}
From \eqref{sec3-eq1}, \eqref{sec3-eq3} and \eqref{sec3-eq4}, $u_n \to u_{\Lambda_s}$ in $L^{\infty}(\mathbb{R}^N)$ as $n \to \infty$. Further, if $u_{\Lambda_s}=0$, then $u_n \to 0$ in $L^{\infty}(\mathbb{R}^N)$ as $n \to \infty$. Moreover, from \ref{f5}, for every $\epsilon>0$ we have
\begin{align*}
    u_n f(u_n) \leq \epsilon u_n^2 \text{ in } \R^N, \text{ for all } n\in \mathbb{N} \text{ large enough}.
\end{align*}
Taking $n\in \mathbb{N}$ large enough, we have
\begin{align*}
    \iint_{\R^N \times \R^N} \frac{|u_n(x)-u_n(y)|^2}{|x-y|^{N+2s}}\,\dxy &= \lambda_n \int_{\R^N} h(x) f(u_n)u_n\,\dx \\
    &\leq \epsilon \lambda_n \int_{\R^N} h(x) u_n^2\,\dx \\
    &\leq \epsilon \frac{\lambda_n}{\lambda_1(h)}\iint_{\R^N \times \R^N} \frac{|u_n(x)-u_n(y)|^2}{|x-y|^{N+2s}}\,\dxy.
\end{align*}
We get $\lambda_n \geq \frac{\lambda_1(h)}{\epsilon}$, for all $n$ large enough. Passing the limit as $n \rightarrow \infty$, we get ${\Lambda_s} \geq \frac{\lambda_1(h)}{\epsilon} \rightarrow +\infty$ when $ \epsilon \rightarrow 0$, a contradiction arises. Thus, $u_{\Lambda_s} \neq 0$ in $\R^N$. Since $u_{\Lambda_s} \in C(\mathbb{R}^N) \cap L^{\infty}(\mathbb{R}^N)$ is nonnegative and satisfies $(-\Delta)^s u_{\Lambda_s} \geq 0$ weakly in $\R^N$. The strong maximum principle \cite[Proposition 5.2.1]{Dipierro-Medina-Valdinoci-2017-Book} infers
$u_{\Lambda_s}>0$ in $\mathbb{R}^N$. 
\end{proof}
%%%%%%%%%%%%%%%%%%%%%%%%%%%%%%%%%%%%%%%%%%%%%%%%%%%%%%%%%%%%
\section{Multiplicity of weak solutions}\label{Section-Multiplicity}
In this section, we prove the existence of at least two distinct positive solutions to \eqref{main-problem} for all $\lambda > \Lambda_s$ using the Classical Linking Theorem. For that purpose, first we need to study the relation between the minimizers of $\Ds2$ and a subspace $X$ of $\Ds2$ defined as follows.
\begin{align}\label{subspace-X}
    X = \left\{ u \in \mathcal{D}^{s,2}(\mathbb{R}^N) \cap C(\mathbb{R}^N) \, : \, \sup_{x \in \mathbb{R}^N} \left( 1 + |x|^{N - 2s} \right) |u(x)|  < \infty \right\}.
\end{align}
The norm on the subspace $X$ is given by
\begin{align*}
    \|u\|_{X} = \sup_{x \in \mathbb{R}^N} \left( 1 + |x|^{N - 2s} \right) |u(x)| .
\end{align*}
\begin{lemma}\label{u-in-subspace}
 Let {\rm\ref{f1}-\ref{f3}}, {\rm\ref{h}, \ref{P1}} and {\rm \ref{P2}} hold.  If $u \in \Ds2$ is a solution to \eqref{main-problem}, then $u \in X$. 
\end{lemma}
\begin{proof}
Let $u \in \Ds2$ be a solution to \eqref{main-problem}. Then $u \in L^\infty(\mathbb{R}^N) \cap C(\mathbb{R}^N)$ using Proposition \ref{regularity}. We fix $g(x):=\lambda h(x) f(u(x))$. The Riesz representation formula for the fractional Laplacian yields
\begin{align*}
u(x) = (-\Delta)^{-s} g(x) = C_{N,-s} \int_{\mathbb{R}^N} \frac{g(y)}{|x - y|^{N - 2s}} \,\dy,
\end{align*}
where $C_{N,-s}$ is a positive constant. The right hand side integral is finite under assumptions \ref{P1}, \eqref{h-bound}, \ref{f4} and using $u \in L^\infty(\R^N)$. For $x \in \mathbb{R}^N \setminus B_1(0)$, we use \eqref{h-bound} and \ref{P2} to estimate $|u(x)|$ as follows:
\begin{align}\label{sec4-eq1}
|u(x)| &\leq C_{N,-s} \int_{\mathbb{R}^N} \frac{|g(y)|}{|x - y|^{N - 2s}} \,\dy \notag\\
&= \lambda C_{N,-s} \int_{\mathbb{R}^N} \frac{h(y)|f(u(y))|}{|x - y|^{N - 2s}} \,\dy \notag\\
&\leq \lambda C_{N,-s} C_0 \int_{\mathbb{R}^N} \frac{P(|y|)|u(y)|}{|x - y|^{N - 2s}} \,\dy \notag\\
&\leq C_1 \lambda \|u\|_{\infty} \int_{\mathbb{R}^N} P(|y|) |x - y|^{2s - N}  \,\dy \notag\\
&\leq C_2 |x|^{2s - N},  \quad \text{ where } C_2= C C_1 \lambda \|u\|_{\infty}.
\end{align}
For $x \in B_1(0)$, we have the following estimate
\begin{align}\label{sec4-eq2}
 |x|^{N - 2s} |u(x)| \leq \|u\|_{\infty}.   
\end{align}
Since $u \in L^\infty(\R^N)$, we conclude from \eqref{sec4-eq1} and \eqref{sec4-eq2} that $u \in X$ i.e., $\sup\limits_{x \in \mathbb{R}^N} \left(1 + |x|^{N - 2s}\right) |u(x)| < \infty$.
\end{proof}
%%%%%%%%%%%%%%%%%%%%%%%%%%%%%%%%%%%%%%%%%%%%%%%%%%
In the following theorem, we establish the relation between the local minimizers of $\mathcal{J}_\lambda$ in $\Ds2$ and those in $X$. A similar result was first proved by Ambrosio~\cite[Theorem~1.1]{Ambrosio-2023-AA} in the fractional setting in $\R^N$ under suitable assumptions on $f$ and $h$. More recently, Carl, Perera, and Tehrani~\cite[Theorem~5.2]{Carl-Perera-Tehrani-2025-Arxiv} extended Ambrosio’s result by considering weaker conditions on $f$ and $h$. Specifically, they assumed that
\begin{align}\label{eq-h}
|h(x)| \leq C_h w(x), \quad \text{where } w(x) = \frac{1}{1 + |x|^{N+\alpha}}, ; x \in \R^N, \text{ for some } \alpha > 0 \text{ and } C_h \geq 0,
\end{align}
and that the function $f : \R \to \R$ is continuous and satisfies the growth condition
\begin{align}\label{eq-f}
|f(t)| \leq C_f \big(1 + |t|^{\gamma-1}\big), \quad 1 \leq \gamma < 2_s^*.
\end{align}
On the one hand our assumptions \ref{f1}--\ref{f5} on $f$ are stronger than condition~\eqref{eq-f}. But on the other hand, the assumptions \ref{h}, \ref{P1}, and \ref{P2} imposed on $h$ are weaker than~\eqref{eq-h}. We have shown a function $h$ in Example \ref{New Example} which is more weaker than \eqref{eq-h}.

\begin{theorem}\label{theorem-4.2}
Assume {\rm \ref{f1}-\ref{f3}, \ref{h}, \ref{P1}} and {\rm \ref{P2}}. If $u\in \Ds2$ is a weak solution to \eqref{main-problem} such that it is a local minimizer of $\mathcal{J}_\lambda$ in $X$, then $u$ is also a local minimizer of $\mathcal{J}_\lambda$ in $\Ds2$.
\end{theorem}
\begin{proof}
Being $u$ a local minimizer of $\mathcal{J}_\lambda$ in the $X$-topology, there exists $\varepsilon_1 > 0$ such that
\begin{align*}
    \mathcal{J}_\lambda(u + v) \geq \mathcal{J}_\lambda(u), \quad \forall v \in X, \, \|v\|_{X} < \varepsilon_1.
\end{align*}
\textbf{Claim:} $u$ is also a local minimizer of $\mathcal{J}_\lambda$ in the $\mathcal{D}^{s,2}(\mathbb{R}^N)$-topology, i.e., there exists $\varepsilon_2 > 0$ such that
\begin{align*}
    \mathcal{J}_\lambda(u + v) \geq \mathcal{J}_\lambda(u), \quad \forall v \in \mathcal{D}^{s,2}(\mathbb{R}^N), \, \|v\|_{\mathcal{D}} < \varepsilon_2.
\end{align*} 
For $n \in \mathbb{N}$, we consider the following minimization problem:
\begin{align*}
 m_n := \inf_{\tilde{u} \in B_n} \mathcal{J}_\lambda(\tilde{u}), \text{ where }   B_n := \left\{ \tilde{u} \in \mathcal{D}^{s,2}(\mathbb{R}^N) : \|\tilde{u} - u\|_{\mathcal{D}} \leq \frac{1}{n} \right\}.
\end{align*}
Notice that $\mathcal{J}_\lambda$ is weakly lower semi-continuous and $B_n$ is weakly closed. Let $\{u_{k,n}\} \subset B_n$ be a minimizing sequence for $\mathcal{J}_\lambda$. Then
\begin{align}\label{sec4-eq3}
 \mathcal{J}_\lambda(u_{k,n}) \to m_n \quad \text{as } k \to \infty. 
\end{align}
By coercivity of $\mathcal{J}_\lambda$ on $\mathcal{D}^{s,2}(\mathbb{R}^N)$ and using \eqref{sec4-eq3}, $\{u_{k,n}\}$ is bounded in $\mathcal{D}^{s,2}(\mathbb{R}^N)$. Therefore, up to a subsequence it converges weakly to some $u_n \in \mathcal{D}^{s,2}(\mathbb{R}^N)$ i.e., $u_{k,n} \rightharpoonup u_n \text{ in } \mathcal{D}^{s,2}(\mathbb{R}^N)$ and $u_n \in B_n$. Thus, we have
\begin{align*}
    m_n \leq \mathcal{J}_\lambda(u_n) \leq \liminf\limits_{k \rightarrow \infty} \mathcal{J}_\lambda(u_{k,n}) = m_n,
\end{align*}
which implies $m_n$ is achieved at some $u_n \in B_n$. There arise two cases: either $\|u_n - u\|_{\mathcal{D}} < \frac{1}{n}$ or $\|u_n - u\|_{\mathcal{D}} = \frac{1}{n}$. If $\|u_n - u\|_{\mathcal{D}} < \frac{1}{n}$, then $u_n$ is a critical point of the functional $\mathcal{J}_\lambda$, i.e.,
\begin{align}\label{sec4-eq04}
    \iint_{\mathbb{R}^N \times \mathbb{R}^N} \frac{(u_n(x) - u_n(y)) (\varphi(x) - \varphi(y))}{|x - y|^{N + 2s}} \, \dxy = \lambda \int_{\mathbb{R}^N} h(x) f(u_n(x)) \varphi(x) \, \dx, \quad \forall \varphi \in \mathcal{D}^{s,2}(\mathbb{R}^N).
\end{align}
If $\|u_n - u\|_{\mathcal{D}} = \frac{1}{n}$, then there exists a Lagrange multiplier $\mu_n \in \mathbb{R}$ such that
\begin{align}\label{sec4-eq4}
    &\iint_{\mathbb{R}^N \times \mathbb{R}^N} \frac{(u_n(x) - u_n(y))(\varphi(x) - \varphi(y))}{|x - y|^{N + 2s}} \, \dxy - \lambda \int_{\mathbb{R}^N} h(x) f(u_n) \varphi(x) \, \dx \notag\\
    &\quad = \mu_n \iint_{\mathbb{R}^N \times \mathbb{R}^N} \frac{((u_n - u)(x) - (u_n - u)(y)) (\varphi(x) - \varphi(y))}{|x - y|^{N + 2s}} \, \dxy,
\quad \text{for all } \varphi \in \mathcal{D}^{s,2}(\mathbb{R}^N).
\end{align}
Following the fact that $u_n$ is a minimizer of $\mathcal{J}_\lambda$ in $B_n$, we have
\begin{align}\label{sec4-eq5}
 \left\langle \mathcal{J}_\lambda'(u_n), u - u_n \right\rangle \geq 0.   
\end{align}
i.e.,
\begin{align*}
    \iint_{\mathbb{R}^N \times \mathbb{R}^N} 
\frac{(u_n(x) - u_n(y)) ((u - u_n)(x) - (u - u_n)(y))}{|x - y|^{N + 2s}} \, \dxy- \lambda \int_{\mathbb{R}^N} h(x) f(u_n) (u - u_n) \, \dx \geq 0.
\end{align*}
Taking $\varphi = u - u_n$ into \eqref{sec4-eq4} and using \eqref{sec4-eq5}, we obtain $\mu_n \leq 0$. Since $u$ is a solution to \eqref{main-problem}, the equation \eqref{sec4-eq4} reduces to
\begin{align*}
    &(1 - \mu_n) \iint_{\mathbb{R}^N \times \mathbb{R}^N} \frac{(u_n(x) - u_n(y))(\varphi(x) - \varphi(y))}{|x - y|^{N + 2s}} \, \dxy 
- \lambda \int_{\mathbb{R}^N} h(x) f(u_n) \varphi(x) \, \dx \\
&\qquad = -\mu_n \iint_{\mathbb{R}^N \times \mathbb{R}^N} 
\frac{(u(x) - u(y)) (\varphi(x) - \varphi(y))}{|x - y|^{N + 2s}} \, \dxy\\
&\qquad = -\lambda \mu_n \int_{\mathbb{R}^N} h(x) f(u) \varphi(x) \, \dx, 
\text{ for all } \varphi \in \mathcal{D}^{s,2}(\mathbb{R}^N).
\end{align*}
Consequently, for every $\varphi \in \Ds2
$ we get
\begin{align}\label{sec4-eq6}
    \iint_{\mathbb{R}^N \times \mathbb{R}^N} \frac{(u_n(x) - u_n(y)) (\varphi(x) - \varphi(y))}{|x - y|^{N + 2s}} \, \dxy 
= \frac{\lambda}{1 - \mu_n} \int_{\mathbb{R}^N} h(x) f(u_n) \varphi(x) \, \dx 
- \frac{\lambda \mu_n}{1 - \mu_n} \int_{\mathbb{R}^N} h(x) f(u) \varphi(x) \, \dx.
\end{align}
Notice that
\begin{align*}
    a_n := \frac{1}{1 - \mu_n} \in (0,1] \text{ and } b_n := \frac{-\mu_n}{1 - \mu_n} \in [0,1).
\end{align*}
For $\|u_n - u\|_{\mathcal{D}} \leq \frac{1}{n}$, we deduce from \eqref{sec4-eq04} and \eqref{sec4-eq6} that $u_n$ is a weak solution to the following problem.
\begin{align}\label{sec4-pde2}
    (-\Delta)^s u_n = \lambda h(x) \left( a_n f(u_n) + b_n f(u) \right)
\quad \text{in } \mathbb{R}^N.
\end{align}
Using $\|u_n - u\|_{\mathcal{D}}\to 0$ as $n \to \infty$ and the Sobolev embedding $\mathcal{D}^{s,2}(\mathbb{R}^N) \hookrightarrow L^{2^*_s}(\mathbb{R}^N)$, we get 
\begin{align}\label{sec4-2star-cnvg}
    u_n \to u \text{ in } L^{2^*_s}(\mathbb{R}^N), \text{ as } n \rightarrow \infty.
\end{align}
There exists a subsequence of $\{u_n\}$, still denoted by itself, and $\zeta \in L^{2^*_s}(\mathbb{R}^N)$ s.t. $|u_n(x)| \leq \zeta(x)$, for a.e. $x \in \mathbb{R}^N$ and for all $n \in \mathbb{N}$. We define
\begin{align}\label{sec4-g(x,w)}
    g_n(x, \omega) = 
\begin{cases}
\lambda h(x) \left( a_n \frac{f(u_n)}{u_n} \omega + b_n f(u) \right), & \text{if } u_n(x) \neq 0, \\
\lambda b_n h(x) f(u), & \text{if } u_n(x) = 0.
\end{cases}
\end{align}
Then $u_n$ is a solution of the problem
\begin{align}\label{sec4-pde3}
 (-\Delta)^s \omega = g_n(x, \omega) \quad \text{in } \mathbb{R}^N .
\end{align}
Since $u \in \mathcal{D}^{s,2}(\mathbb{R}^N)$ is a weak solution of \eqref{main-problem}, we get $u \in X$ due to Lemma \ref{u-in-subspace}. Further, we have
\begin{align*}
    |f(u(x))| \leq C_0 \|u\|_\infty, \text{ for all } x \in \mathbb{R}^N.
\end{align*}
From \ref{f4}, we observe that
\begin{align*}
    \frac{f(u_n)}{u_n} \leq C_0, \quad \forall n \in \mathbb{N} \text{ and for all } x \in \{ u_n \neq 0 \}.
\end{align*}
Thus, we have
\begin{align*}
   | g_n(x, \omega)| \leq C_0 \lambda h(x) (|\omega| + \|u\|_\infty), \quad \text{ for all } x \in \mathbb{R}^N.
\end{align*}
Since $u_n$ solves weakly to \eqref{sec4-pde3}, we follow arguments similar to Proposition \ref{regularity} and obtain
\begin{align*}
    \|u_n\|_{L^\tau(\mathbb{R}^N)} \leq C_1 + C_2 \|u_n\|_{L^{2^*_s}(\mathbb{R}^N)}, \quad \text{ for all } n \in \mathbb{N},
\end{align*}
where $\tau \geq 2_s^*$ and $C_1, C_2 > 0$. From \eqref{sec4-2star-cnvg}, $\{u_n\}$ is bounded in $L^\tau(\mathbb{R}^N)$ for all $\tau \geq 2^*_s$. Applying the Moser iteration as in Proposition \ref{regularity}, we deduce that $u_n \in C(\mathbb{R}^N) \cap L^\infty(\mathbb{R}^N)$.

Let $w_n = u_n - u$. Since $u$ is a solution to \eqref{main-problem}, we get
\begin{align*}
(-\Delta)^s w_n &= (-\Delta)^s u_n - (-\Delta)^s u \\
&= \lambda h(x) \left( a_n f(u_n) + b_n f(u) \right) - \lambda h(x) f(u) \\
&= \lambda h(x) a_n f(u_n) + \left(b_n - 1 \right) \lambda h(x) f(u) \\
&= \lambda a_n h(x) \left( f(u_n) - f(u) \right) \\
&= \lambda a_n h(x) \, G(x, u_n) \, \omega_n \quad \text{in } \mathbb{R}^N,
\end{align*}
where
\begin{align*}
    G(x, u_n) = 
\begin{cases}
\frac{f(u_n) - f(u)}{u_n - u}, & \text{if } u_n \neq u, \\
0, & \text{if } u_n = u.
\end{cases}
\end{align*}
Since $f$ is locally Lipschitz in $\R$ (see \ref{f3}), we get $|G(x, w_n)| \leq L$, for sufficiently large $n \in \mathbb{N}$ and for some $L>0$. We further define $\widetilde{G}(x, w_n) := \lambda a_n h(x) G(x, w_n) w_n$. Using \( a_n \in (0,1] \), we find that $|\widetilde{G}(x, w_n)| \leq \lambda L h(x) |w_n |$. Following the arguments as in Proposition \ref{regularity}, we get
\begin{align}\label{sec4-wn-estimate-L-infty}
    \|w_n\|_{\infty} \leq C \|w_n \|_{2^*_s}.
\end{align}
The Riesz potential formula for the fractional Laplacian yields
\begin{align}\label{sec4-Riesz-1}
    |w_n(x)| &= \left| (-\Delta)^{-s}(\lambda a_n h(x) G(x, w_n) w_n) \right| \leq C_{N,-s} \, \lambda L \int_{\mathbb{R}^N} \frac{h(y) |w_n(y)|}{|x - y|^{N - 2s}} \,\dy \notag\\
    & \leq C_1 \|w_n\|_{\infty} \int_{\mathbb{R}^N} \frac{P(|y|)}{|x - y|^{N - 2s}} \,\dy \notag\\
    & \leq C_2 \|w_n\|_{\infty} |x|^{2s - N}, \quad \text{ for all } x \in \mathbb{R}^N \setminus B_1(0).
\end{align}
For $x \in B_1(0)$, we have the following estimate
\begin{align}\label{sec4-Riesz-2}
    |x|^{N - 2s} |w_n(x)| \leq \|w_n\|_{\infty}.
\end{align}
Combining the estimates \eqref{sec4-Riesz-1}, \eqref{sec4-Riesz-2} and the fact that $w_n \in L^\infty(\R^N)$, we deduce that
\begin{align}\label{sec4-eq7}
    \left(1 + |x|^{N - 2s} \right) |w_n(x)| \leq C_3 \|w_n\|_{\infty} , \text{ for some } C_3>0.
\end{align}
Using \eqref{sec4-2star-cnvg} and \eqref{sec4-wn-estimate-L-infty}, $\|w_n\|_\infty \rightarrow 0$ as $n \rightarrow \infty$ and consequently we get $\|w_n\|_X \rightarrow 0$ as $n \rightarrow \infty$ due to the estimate \eqref{sec4-eq7}. Hence, for $n$ large enough, we get
\begin{align*}
    \inf_{ \tilde{u} \in B_n} \mathcal{J}_\lambda( \tilde{u}) = \mathcal{J}_\lambda(u_n) = \mathcal{J}_\lambda\left( u + (u_n - u) \right) \geq \mathcal{J}_\lambda(u),
\end{align*}
i.e., $u$ is a local minimizer of $\mathcal{J}_\lambda$ in $\mathcal{D}^{s,2}(\mathbb{R}^N)$-topology.
\end{proof}
\begin{lemma}\label{Lemma-4.3}
  The function $u \equiv 0$ is a local minimizer of $\mathcal{J}_\lambda$ in $(X, \|\cdot\|_X)$.  
\end{lemma}
\begin{proof}
Using \ref{f5}, for each $\lambda > 0$, there exists $t_\lambda > 0$ such that
\begin{align*}
    f(t) \leq \frac{\lambda_1(h)}{2\lambda } t, \quad t \in (0, t_\lambda).
\end{align*}
Then
\begin{align*}
    F(t) = \int_0^t f(\tau)\, \mathrm{d}\tau \leq \frac{\lambda_1(h)}{2\lambda } \int_0^t \tau\, \mathrm{d}\tau = \frac{\lambda_1(h)}{4\lambda } t^2, \quad t \in (0, t_\lambda)
\end{align*}
and
\begin{align*}
    \int_{\mathbb{R}^N} h(x) F(u)\, \dx \leq \frac{\lambda_1(h)}{4\lambda } \int_{\mathbb{R}^N} h(x) u^2\, \dx,
\end{align*}
for all $u \in X$ such that $|u|_\infty \leq \|u\|_X \leq t_\lambda$. Finally,
\begin{align*}
    \begin{aligned}
\mathcal{J}_\lambda(u) &= \frac{1}{2} \iint_{\mathbb{R}^N \times \mathbb{R}^N} \frac{|u(x) - u(y)|^2}{|x - y|^{N + 2s}} \, \dxy - \lambda \int_{\mathbb{R}^N} h(x) F(u)\, \dx \\
&\geq \frac{1}{2} \iint_{\mathbb{R}^N \times \mathbb{R}^N} \frac{|u(x) - u(y)|^2}{|x - y|^{N + 2s}} \, \dxy - \frac{\lambda_1(h)}{4 } \int_{\mathbb{R}^N} h(x) u^2\, \dx \\
&\geq \frac{1}{4} \iint_{\mathbb{R}^N \times \mathbb{R}^N} \frac{|u(x) - u(y)|^2}{|x - y|^{N + 2s}} \, \dxy > 0,
\end{aligned}
\end{align*}
for all $u \in B_{t_\lambda}^X(0) \setminus \{0\}$, where $B_{t_\lambda}^X(0)$ is the ball in $X$ with radius $t_\lambda$ and center at origin. 
\end{proof}
Now our aim is to find a second local minimizer of $\mathcal{J}_\lambda$ with respect to $(X,\|\cdot\|_X)$ for $\lambda > {\Lambda_s}$. Let $\lambda_2>\lambda>\lambda_1>\Lambda_s$ and $u_{\lambda_1}, u_{\lambda_2}$ be the weak solutions of $(P_{\lambda_1}^s)$ and $(P_{\lambda_2}^s)$, respectively. The solutions $u_{\lambda_1}$ and $ u_{\lambda_2}$ are positive by Proposition \ref{Positivity}. Moreover, following similar arguments used to get \eqref{solution-compare}, we obtain
\begin{align*}
    0 < u_{\lambda_1} \leq u_{\lambda_2} \quad \text{in } \mathbb{R}^N.
\end{align*}
The monotonicity of $f$ due to \ref{f3} yields
\begin{align*}
    0 < f(u_{\lambda_1}) \leq f(u_{\lambda_2}) \quad \text{in } \mathbb{R}^N.
\end{align*}
Let us denote $w = u_{\lambda_2} - u_{\lambda_1}$. Then, $w$ satisfies the following inequality in the weak sense.
\begin{align*}
(- \Delta)^s w &= (- \Delta)^s u_{\lambda_2}  - (- \Delta)^s u_{\lambda_1} \\
&= \lambda_2 h(x) f(u_{\lambda_2}) - \lambda_1 h(x) f(u_{\lambda_1})\\
&> \lambda_1 h(x)\left(f(u_{\lambda_2}) - f(u_{\lambda_1})\right)\\
&\geq 0 \quad \text{in } \mathbb{R}^N.
\end{align*}
Since $w \geq 0$ in $\mathbb{R}^N$ and $w \not\equiv 0$ in $\mathbb{R}^N$, the strong maximum principle infers $w > 0$ in $\mathbb{R}^N$ i.e., $u_{\lambda_2} > u_{\lambda_1}$ in $\mathbb{R}^N$. Further, we define a function $\widehat{f}: \R \to \R$ as follows:
\begin{align*}
    \widehat{f}(t) = 
\begin{cases}
f(u_{\lambda_1}) & \text{if } t \leq u_{\lambda_1}, \\
f(t) & \text{if } u_{\lambda_1} \leq t \leq u_{\lambda_2}, \\
f(u_{\lambda_2}) & \text{if } t \geq u_{\lambda_2},
\end{cases}
\end{align*}
and its primitive is given by $\widehat{F}(t) = \int_0^t \widehat{f}(\tau)\, \mathrm{d}\tau$. The associated energy functional $\widehat{\mathcal{J}}_\lambda: \Ds2
 \to \R$ is given by
\begin{align*}
    \widehat{\mathcal{J}}_\lambda(u) = \frac{1}{2} \iint_{\mathbb{R}^N \times \mathbb{R}^N} \frac{|u(x)-u(y)|^2}{|x-y|^{N+2s}}\,\dxy - \lambda \int_{\mathbb{R}^N} h(x) \widehat{F}(u)\, \dx.
\end{align*}
Now $\widehat{f}$ is a bounded and continuous in $\mathbb{R}$ which implies that $\widehat{\mathcal{J}}_\lambda$ admits a global minimizer $v_\lambda$ in $\mathcal{D}^{s,2}(\mathbb{R}^N)$. As a result, $v_\lambda$ weakly solves the problem
\begin{equation}\label{sec4-P-lambda-hat-problem}
\left\{
\begin{aligned}
(- \Delta)^s u &= \lambda h(x) \widehat{f}(u) \quad \text{in } \mathbb{R}^N, \\
u &> 0 \quad \text{in } \mathbb{R}^N, \quad u \in \mathcal{D}^{s,2}(\mathbb{R}^N).
\end{aligned}
\right.
\end{equation}
From Proposition \ref{regularity}, we get $v_\lambda \in C(\mathbb{R}^N) \cap L^\infty(\mathbb{R}^N)$. Next our goal is show that $v_\lambda \in [u_{\lambda_1}, u_{\lambda_2}]$. Observe that
\begin{align*}
    \iint_{\mathbb{R}^N \times \mathbb{R}^N} \frac{(u_{\lambda_2}(x)-u_{\lambda_2}(y))(\varphi(x)- \varphi(y))}{|x-y|^{N+2s}}\,\dxy  = \lambda_2 \int_{\mathbb{R}^N} h(x) f(u_{\lambda_2}) \varphi\, \dx \geq \lambda \int_{\mathbb{R}^N} h(x) f(u_{\lambda_2}) \varphi\, \dx
\end{align*}
and
\begin{align*}
    \iint_{\mathbb{R}^N \times \mathbb{R}^N} \frac{(v_{\lambda}(x)-v_{\lambda}(y))(\varphi(x)- \varphi(y))}{|x-y|^{N+2s}}\,\dxy  = \lambda \int_{\mathbb{R}^N} h(x) \widehat{f}(v_\lambda) \varphi\, \dx,
\end{align*}
for all $\varphi \in \mathcal{D}^{s,2}(\mathbb{R}^N)$ and $\varphi \geq 0$. Therefore, we get
\begin{align*}
    \iint_{\mathbb{R}^N \times \mathbb{R}^N} \frac{\left( \left(v_{\lambda}(x) -u_{\lambda_2}(x) \right)-\left(v_{\lambda}(y) -u_{\lambda_2}(y) \right)\right)(\varphi(x)- \varphi(y))}{|x-y|^{N+2s}}\,\dxy \leq \lambda \int_{\mathbb{R}^N} h(x)\left(\widehat{f}(v_\lambda) - f(u_{\lambda_2})\right) \varphi\, \dx.
\end{align*}
Taking $\varphi = (v_\lambda - u_{\lambda_2})^+ = \max\{v_\lambda - u_{\lambda_2}, 0\} \geq 0$, we have
\begin{align*}
   \iint_{\mathbb{R}^N \times \mathbb{R}^N} \frac{\left( (v_\lambda - u_{\lambda_2})^+(x)-(v_\lambda - u_{\lambda_2})^+(y)\right)^2}{|x-y|^{N+2s}}\,\dxy \leq \lambda \int_{\mathbb{R}^N} h(x)\left(\widehat{f}(v_\lambda) - f(u_{\lambda_2})\right)(v_\lambda - u_{\lambda_2})^+\, \dx = 0. 
\end{align*}
This implies $\|(v_\lambda - u_{\lambda_2})^{+}\|_{\mathcal{D}} \leq 0$, which further gives
$(v_\lambda - u_{\lambda_2})^+ = 0$, i.e., $v_\lambda \leq u_{\lambda_2}$ in $\mathbb{R}^N$. Repeating the same argument as before, we have $v_\lambda < u_{\lambda_2}$ in $\mathbb{R}^N$. Similarly, we can show that $u_{\lambda_1} < v_\lambda$ in $\mathbb{R}^N$, which gives $\widehat{f}(v_\lambda) = f(v_\lambda)$ and consequently $v_\lambda$ is also a solution of \eqref{main-problem}. Following the similar arguments as given in Lemma \ref{u-in-subspace}, we immediately conclude that $u_{\lambda_1}, v_\lambda, u_{\lambda_2} \in X$. 

%For the sake of simplicity, we will prove that $v_\lambda \in X$ and the proof for the other two functions follows using similar arguments. The Riesz Potential Formula for the fractional Laplacian gives
% \begin{align*}
%     v_\lambda(x) = \lambda C_{N,s} \int_{\mathbb{R}^N} \frac{h(y) f(v_\lambda(y))}{|x - y|^{N - 2s}}\, \dy, \quad x \in \mathbb{R}^N, \text{ for some } C_{N,s} > 0.
% \end{align*}
% From \ref{P2}, we establish the following estimate:
% \begin{align*}
% |v_\lambda(x)| &\leq C(N,s, \lambda) \|u_{\lambda_2}\|_\infty \int_{\mathbb{R}^N} \frac{P(|y|)}{|x - y|^{N - 2s}}\, \dy \notag \\
% &\leq C(N,s, \lambda, \|u_{\lambda_2}\|_\infty) |x|^{2s - N}, \quad \text{ for all }|x| \geq 1. 
% \end{align*}
% Since $v_\lambda \in L^\infty(\R^N)$, we have $|x|^{N - 2s}|v_\lambda(x)| \leq \|v_\lambda\|_\infty$, for all $|x|<1$. Subsequently, for every $x \in \R^N$ we get
% \begin{align*}
%     (1 + |x|^{N - 2s})|v_\lambda(x)| \leq 2 \|v_\lambda\|_\infty + C(N,s,\lambda, \|u_{\lambda_2}\|_\infty) < \infty,
% \end{align*}
% which implies $v_\lambda \in X$. 
%%%%%%%%%%%%%%%%%%%%%%%%%%%%%%%%%%%%%%%%%%%%%%%%%%%%%%%%%%%%
For $p\in (1, \frac{N}{s})$ and an open set \( \Omega \subset \mathbb{R}^N \), we define
\begin{align*}
    \widetilde{\mathcal{D}}^{s,p}(\Omega) := \left\{ u \in L^{p-1}_{loc}(\mathbb{R}^N) \cap L^{p^*_{s}}(\Omega) \, : \, \exists E \supset \Omega \text{ with } E^c \text{ compact, } \text{dist}(E^c, \Omega) > 0 \text{ and } [u]^p_{W^{s,p}(E)} < +\infty \right\},
\end{align*}
where
    \begin{align*}
    [u]^p_{W^{s,p}(E)} := \iint_{E \times E} \frac{|u(x) - u(y)|^p}{|x - y|^{N + sp}}\,\dxy.
    \end{align*}  
Notice that $\|u\|_{\mathcal{D}} = [u]_{W^{s,p}(\mathbb{R}^N)}$.    
\begin{lemma}{\rm(cf. \cite[Proposition 2.5]{Brasco-optimal-decay})}\label{Proposition-2.5}
    For any \( u \in \widetilde{\mathcal{D}}^{s,p}(\Omega) \), the operator $(-\Delta_p)^su$ given by
    \begin{align*}
    \mathcal{D}_0^{s,p}(\Omega) \ni \varphi \mapsto \left\langle (-\Delta_p)^s u, \varphi \right\rangle := \iint_{\mathbb{R}^N \times \mathbb{R}^N} \frac{\mathcal{J}_p(u(x) - u(y))(\varphi(x) - \varphi(y))}{|x - y|^{N + sp}}\,\dxy
    \end{align*}
    is well-defined and belongs to the dual space \( \left( \mathcal{D}_0^{s,p}(\Omega) \right)^* \), where \( \mathcal{J}_p(t) = |t|^{p-2}t \), for $t \in \R$. The space $\mathcal{D}_0^{s,p}(\Omega)$ is defined as follows (cf. \cite[Theorem 2.1]{Brasco-optimal-decay})
    \begin{align*}
        \mathcal{D}_0^{s,p}(\Omega):=\left\{u \in L^{p_s^*}(\Omega) : u \equiv 0 \text{ in } \Omega^c , \|u\|_{\mathcal{D}}< + \infty\right\}.
    \end{align*}
    Moreover, if $\partial \Omega$ is compact and locally the graph of a continuous function, then $\mathcal{D}_0^{s,p}(\Omega)$ is the completion of $C_c^\infty(\Omega)$ with respect to the norm $\|\cdot\|_{\mathcal{D}}$.
\end{lemma}   
\begin{definition}{\rm(cf. \cite[Definition 2.6]{Brasco-optimal-decay})}\label{Definition-2.6}
    Let \(  u \in \widetilde{\mathcal{D}}^{s,p}(\Omega)  \) and \( \Lambda \in \left( \mathcal{D}_0^{s,p}(\Omega) \right)^* \). We say that $(-\Delta_p)^s u \leq \Lambda$ weakly in $\Omega$ if for all \( \varphi \in \mathcal{D}_0^{s,p}(\Omega), \varphi \geq 0 \text{ in } \Omega \),
    \begin{align*}
    \iint_{\mathbb{R}^N \times \mathbb{R}^N} \frac{\mathcal{J}_p(u(x) - u(y))(\varphi(x) - \varphi(y))}{|x - y|^{N + sp}}\,\dxy \leq \langle \Lambda, \varphi \rangle.
    \end{align*}
\end{definition}
\begin{lemma}{\rm(cf. \cite[Theorem A.4]{Brasco-optimal-decay})}\label{Theorem-A.4}
    For any \( R > 0 \) and $p \in (1, \frac{N}{s})$ \( \Gamma_p(x) = |x|^{-\frac{N - sp}{p - 1}} \) belongs to \(\widetilde{\mathcal{D}}^{s,p}(\overline{B}_R^c) \) and weakly solves
     \begin{align*}
          (-\Delta_p)^s u = 0 \quad \text{in } \overline{B}_R^c.
     \end{align*}
\end{lemma}
%%%%%%%%%%%%%%%%%%%%%%%%%%%%%%%%%%%%%%%%%%%%%%%%%%%%%%
\begin{proposition}\label{Prop-4.7}
    The function $v_\lambda$ is a local minimizer of $\mathcal{J}_\lambda$ in $X$.
\end{proposition}
\begin{proof}
Since $u_{\lambda_2}$ and $v_\lambda$ are the solutions to $(P_{\lambda_2}^s)$ and \eqref{main-problem} respectively, we get $u_{\lambda_2}, v_\lambda \in X$. There exist $\delta>0$ and $R_\delta>0$ such that $ |u_{\lambda_2}(x)| < \delta$, $|v_\lambda(x)|<\delta$ whenever $|x| \geq R_\delta$. We define $z = u_{\lambda_2} - v_\lambda$, which implies $z>0$ in $\R^N$. Taking $\mu = \min\{z(x): x \in \overline{B(0,R_\delta)} \}>0$, we have
\begin{align*}
    z(x) \geq \mu > \frac{\mu}{2(1+|x|^{N-2s})}, \quad \text{ for all } x \in \overline{B(0,R_\delta)}.
\end{align*}     
Note that $(-\Delta)^s z \geq 0$ weakly in $\R^N$ i.e.
\begin{align}\label{sec4-w1}
    \iint_{\mathbb{R}^N \times \mathbb{R}^N} \frac{(z(x) - z(y)) (\varphi(x) - \varphi(y))}{|x - y|^{N+2s}} \, \dxy \geq 0,
    \quad \text{for every } \varphi \in \mathcal{D}^{s,p}(\mathbb{R}^N) \text{ and } \varphi \geq 0.
\end{align}
Pick $\nu>0$ such that $z(x) >\Gamma_2(x)$ for $|x| = R_\delta$, where $ \Gamma_2(x)=\nu |x|^{-(N-2s)}$ . Let $\Omega_1:= \R^N \setminus \overline{B(0,R_\delta)}$. From Lemma \ref{Theorem-A.4}, we get $\Gamma_2 \in \widetilde{\mathcal{D}}^{s,2}(\Omega_1)$ (cf. Definition \ref{Definition-2.6}) and $\Gamma_2$ weakly solves $(-\Delta)^s u =0$ in $\Omega_1$ i.e.,
\begin{align}\label{sec4-w2}
    \iint_{\mathbb{R}^N \times \mathbb{R}^N} \frac{(\Gamma_2(x) - \Gamma_2(y))(\varphi(x) - \varphi(y))}{|x - y|^{N + 2s}} \, \dxy = 0, \text{ for all } \varphi \in \mathcal{D}_0^{s,2}(\Omega_1).
\end{align}
From \eqref{sec4-w1} and \eqref{sec4-w2}, we deduce that
\begin{align}\label{sec4-w3}
    \iint_{\mathbb{R}^N \times \mathbb{R}^N}
    \frac{\left( z(x) - z(y) - \left((\Gamma_2(x) - \Gamma_2(y)\right) \right) (\varphi(x) - \varphi(y))}{|x - y|^{N + 2s}} \, \dxy \geq 0, \text{ for all }
  \varphi \in \mathcal{D}_0^{s,2}(\Omega_1), \; \varphi \geq 0.
\end{align}
Let $\varphi = \left(z - \Gamma_2 \right)^- := \max\left\{ 0 , - \left( z - \Gamma_2 \right)\right\}$ in $\Omega_1$, which can be extended by zero in $\overline{B(0,R_\delta)}$ using \cite[Theorem 5.4]{Hitchiker-guide-2012}. Let us denote $\Psi= \left(z - \Gamma_2 \right)^+ :=\max\left\{ 0 ,  z - \Gamma_2 \right\}$. Then
\begin{align}\label{sec4-w4}
  &\left( z(x) - z(y) - \left((\Gamma_2(x) - \Gamma_2(y)\right) \right) (\varphi(x) - \varphi(y)) \notag\\  &=\left(\Psi(x) - \Psi(y) - (\varphi(x)- \varphi(y) \right) \left( \varphi(x) -\varphi(y) \right) \notag\\
  &= (\Psi(x) - \Psi(y))\left( \varphi(x) -\varphi(y) \right) - \left( \varphi(x) -\varphi(y) \right)^2 \notag\\
    &= \Psi(x)\varphi(x) - \Psi(y)\varphi(x) - \Psi(x)\varphi(y) + \Psi(y)\varphi(y)- \left( \varphi(x) -\varphi(y) \right)^2 \notag\\
    &= - \Psi(y)\varphi(x) - \Psi(x)\varphi(y) - \left( \varphi(x) -\varphi(y) \right)^2 \notag\\
    & \leq - \left( \varphi(x) -\varphi(y) \right)^2.
\end{align}
From \eqref{sec4-w3} and \eqref{sec4-w4}, we deduce that
\begin{align}\label{sec4-w5}
    -\iint_{\R^N \times \R^N} \frac{\left( \varphi(x) -\varphi(y) \right)^2}{|x-y|^{N+2s}}\,\dxy \geq 0.
\end{align}
We estimate the above integral as follows
\begin{align}\label{sec4-w6}
   & \iint_{\mathbb{R}^N \times \mathbb{R}^N} \frac{(\varphi(x) - \varphi(y))^2}{|x - y|^{N + 2s}} \,\dxy  \notag \\
   &= \iint_{\Omega_1 \times \Omega_1} \frac{(\varphi(x) - \varphi(y))^2}{|x - y|^{N + 2s}} \,\dxy + 2 \iint_{\Omega_1^c \times \Omega_1} \frac{(\varphi(y))^2}{|x - y|^{N + 2s}} \,\dxy + \iint_{\Omega_1^c \times \Omega_1^c} \frac{(\varphi(x) - \varphi(y))^2}{|x - y|^{N + 2s}} \,\dxy \notag\\
    & \geq \iint_{\Omega_1 \times \Omega_1} \frac{(\varphi(x) - \varphi(y))^2}{|x - y|^{N + 2s}} \dxy.
\end{align}
The estimates \eqref{sec4-w6} and \eqref{sec4-w5} yield $\varphi = 0$ in $\Omega_1$. Subsequently, we get
\begin{align}
    z(x) \geq \frac{\nu}{|x|^{N - 2s}} > \frac{\nu}{2(1+|x|^{N - 2s})} \quad \text{in } \Omega_1 = \R^N \setminus \overline{B(0,R_\delta)}.
\end{align}
Choosing $\epsilon_0 = \min\left\{\frac{\mu}{2}, \frac{\nu}{2} \right\}$ such that $ z(x) > \frac{\epsilon_0}{1+|x|^{N - 2s}}$ for all $x \in \R^N$ and which further gives
\begin{align}\label{in1}
    v_\lambda(x) < v_\lambda(x) + \frac{\epsilon_0}{1+|x|^{N - 2s}} < u_{\lambda_2}(x), \text{ for all } x\in  \R^N.
\end{align}
Using similar set of arguments, we also get
\begin{align}\label{in2}
    u_{\lambda_1}(x) < v_\lambda(x) - \frac{\epsilon_0}{1+|x|^{N - 2s}} < v_{\lambda}(x), \text{ for all } x\in  \R^N.
\end{align}
Let $B_{\epsilon_0}^X(v_\lambda)$ denote an open ball in $X$ centered at $v_\lambda$ and with radius $\epsilon_0$, i.e., $B_{\epsilon_0}^X(v_\lambda)= \left\{ v \in X : \|v-v_\lambda\|_X<\epsilon_0 \right\}$. Let $u \in B_{\epsilon_0}^X(v_\lambda)$. Then $u$ satisfies the following inequality
\begin{align*}\label{in3}
     v_\lambda(x) - \frac{\epsilon_0}{1+|x|^{N - 2s}}  < u(x)< v_\lambda(x) + \frac{\epsilon_0}{1+|x|^{N - 2s}}, \text{ for all } x\in  \R^N.
\end{align*}
We combine the above inequality with the inequalities 
\eqref{in1} and \eqref{in2}, and deduce the following
\begin{align}
   u_{\lambda_1}(x) < u(x)< u_{\lambda_2}(x), \text{ for all } x\in  \R^N.
\end{align}
Let us denote $K = \int_{\R^N}h(x)\left(f(u_{\lambda_1})u_{\lambda_1}- F(u_{\lambda_1}) \right)\,\dx$. Next we claim that $\mathcal{J}_\lambda(u) = \Hat{\mathcal{J}_\lambda}(u)+ \lambda K$, for every $u \in B_{\epsilon_0}^X(v_\lambda)$. Observe that
\begin{align*}
    \Hat{\mathcal{J}_\lambda}(u)+ \lambda K &= \Hat{\mathcal{J}_\lambda}(u)+ \lambda \int_{\R^N}h(x)\left(f(u_{\lambda_1})u_{\lambda_1}- F(u_{\lambda_1}) \right)\,\dx \\
    &= \frac{1}{2}\|u\|_{\mathcal{D}}^2 - \lambda \int_{\R^N}h(x) \Hat{F}(u)\,\dx + \lambda \int_{\R^N}h(x)\left(f(u_{\lambda_1})u_{\lambda_1}- F(u_{\lambda_1}) \right)\,\dx \\
    &= \frac{1}{2}\|u\|_{\mathcal{D}}^2 - \lambda \int_{\R^N}h(x) \left( \int_{0}^{u_{\lambda_1}}\Hat{f}(s)\,\ds + \int_{u_{\lambda_1}}^{u}\Hat{f}(s)\,\ds \right)\,\dx + \lambda \int_{\R^N}h(x)\left(f(u_{\lambda_1})u_{\lambda_1}- F(u_{\lambda_1}) \right)\,\dx \\
    &=\frac{1}{2}\|u\|_{\mathcal{D}}^2 - \lambda \int_{\R^N}h(x) F(u)\,\dx\\
    &= \mathcal{J}_\lambda(u). 
\end{align*}
Since $v_\lambda$ is a global minimizer of $\Hat{\mathcal{J}_\lambda}$ in $\Ds2$, we have
\begin{align*}
    \mathcal{J}_\lambda(u) = \Hat{\mathcal{J}_\lambda}(u)+ \lambda K \geq \Hat{\mathcal{J}_\lambda}(v_\lambda)+ \lambda K = \mathcal{J}_\lambda(v_\lambda), \text{ for all } u \in B_{\epsilon_0}^X(v_\lambda).
\end{align*}
This completes the proof of this Lemma.
\end{proof}
\begin{remark}\label{Rem-2}
From Theorem \ref{theorem-4.2}, Lemma \ref{Lemma-4.3} and Proposition \ref{Prop-4.7}, we conclude that $0$ and $v_\lambda$ are local minimizers of $\mathcal{J}_\lambda$ in $\Ds2$.
\end{remark}
We know that $u_\lambda$ and $v_\lambda$ are solutions to \eqref{main-problem}, but we do not have any information about whether these solutions are distinct or not. We just know that $v_\lambda \in (u_{\lambda_1}, u_{\lambda_2})$. In the following, we construct a distinct solution to \eqref{main-problem} using the first solution $u_\lambda$. We define a functional $\mathcal{I}_\lambda : \Ds2 \rightarrow \R$ given by{\small  
\begin{align*}
    \mathcal{I}_\lambda (w) &= \frac{1}{2}\iint\limits_{\R^{2N}} \frac{|(v_\lambda +w)(x)-(v_\lambda + w)(y)|^2}{|x-y|^{N+2s}}\,\dxy - \iint\limits_{\R^{2N}} \frac{(v_\lambda(x)-v_\lambda(y)) (w(x)-w(y))}{|x-y|^{N+2s}}\,\dxy- \lambda \int\limits_{\R^N} h(x) G(x,w)\,\dx,
\end{align*}}
where 
\begin{align*}
G(x,t)= \int_{0}^{t}g(x,\sigma)\,\mathrm{d}\sigma, \quad \text{ and } \quad   g(x,\sigma) = \begin{cases}
        f(v_\lambda (x) +\sigma) - f(v_\lambda (x) ) , & \sigma<0,\\
        0,& \sigma \geq 0.
    \end{cases}
\end{align*}
The functional $\mathcal{I}_\lambda$ is Frech\'{e}t differentiable on $\Ds2$ and its Frech\'{e}t derivative is given by
\begin{align*}
    \langle \mathcal{I}_{\lambda}'(w), \varphi \rangle &=\iint_{\R^N \times \R^N} \frac{\left((v_\lambda +w)(x)-(v_\lambda + w)(y)\right) (\varphi(x)-\varphi(y))}{|x-y|^{N+2s}}\,\dxy - \lambda \int_{\R^N} h(x) g(x,w)\varphi\,\dx\\
    &\qquad -  \iint_{\R^N \times \R^N} \frac{(v_\lambda(x)-v_\lambda(y)) (\varphi(x)-\varphi(y))}{|x-y|^{N+2s}}\,\dxy,\quad \text{ for all } \varphi \in \Ds2.
\end{align*}
\begin{lemma}
    $w_1= 0$ and $w_2 = -v_\lambda$ are local minimizers of $\mathcal{I}_\lambda$ in $\Ds2$.
\end{lemma}
\begin{proof}
For each $x\in \R^N$ and $t\in \R$ such that $t<0$, we get
    \begin{align*}
        G(x,t)&= \int_{0}^{t} f(v_\lambda (x) +\sig)\,\mathrm{d}\sigma - \int_{0}^{t} f(v_\lambda (x))\,\mathrm{d}\sigma \\
        &= \int_{v_\lambda}^{v_\lambda + t} f(r)\,\mathrm{d}r - f(v_\lambda(x))t \qquad  \text{(Putting $r= v_\lambda +\sigma$)}\\
        &= F(v_\lambda (x) + t) - F(v_\lambda (x)) -f(v_\lambda (x))t.
    \end{align*}
For $t\in \R$, we recall that $t=t^+ + t^-$, where $t^+= \max\{t,0\}$ and $t^- = \min\{t,0\}$. Notice that $G(x, t^+)=0$, for all $t \in \R$ and $x\in \R^N$. For any $x \in \R^N$ and $t\in \R$, we get
\begin{align}\label{sec4-eq8}
    G(x,t) = F(v_\lambda (x) + t^-) - F(v_\lambda (x)) -f(v_\lambda (x))t^-.
\end{align}
Let $w \in \Ds2$ and $\delta>0$ be sufficiently small so that $\|w\|_{\mathcal{D}}< \delta$. Then
\begin{align}\label{sec4-eq011}
    \mathcal{I}_\lambda (w) - \mathcal{I}_\lambda (0) &= \frac{1}{2} \iint_{\R^N \times \R^N} \frac{|(v_\lambda +w) (x) - (v_\lambda +w) (y) |^2}{|x-y|^{N+2s}}\,\dxy - \iint_{\R^N \times \R^N} \frac{(v_\lambda (x) - v_\lambda (y)) (w(x) - w(y))}{|x-y|^{N+2s}}\,\dxy \notag\\
    &\quad - \lambda \int_{\R^N} h(x)G(x,w)\,\dx - \frac{1}{2} \iint_{\R^N \times \R^N} \frac{|v_\lambda (x) - v_\lambda  (y) |^2}{|x-y|^{N+2s}}\,\dxy.
\end{align}
Notice that $w= w^+ + w^-$, with $w^+ \geq 0 $ and $ w^- \leq 0$. We estimate the following expression:
\begin{align} \label{sec4-eq11}
    B &= \left( (v_\lambda +w) (x) - (v_\lambda +w) (y) \right)^2 \notag\\
    &= \left( (v_\lambda +w^-) (x) - (v_\lambda +w^-) (y) + \left(w^+(x) - w^+(y) \right) \right)^2 \notag\notag\\
    &= \left((v_\lambda +w^-) (x) - (v_\lambda +w^-) (y) \right)^2 + \left( w^+(x) - w^+(y) \right)^2 + 2 \left((v_\lambda +w^-) (x) - (v_\lambda +w^-) (y) \right)\left( w^+(x) - w^+(y) \right) \notag\\
    & = \left((v_\lambda +w^-) (x) - (v_\lambda +w^-) (y) \right)^2 + \left( w^+(x) - w^+(y) \right)^2 + 2 \left(v_\lambda (x) - v_\lambda (y) \right)\left( w^+(x) - w^+(y) \right) \notag\\
    & \qquad - 2 w^-(x)w^+(y) - 2 w^-(y) w^+(x) \notag\\
    & \geq \left((v_\lambda +w^-) (x) - (v_\lambda +w^-) (y) \right)^2 + \left( w^+(x) - w^+(y) \right)^2 + 2 \left(v_\lambda (x) - v_\lambda (y) \right)\left( w^+(x) - w^+(y) \right) \\
    & \geq \left((v_\lambda +w^-) (x) - (v_\lambda +w^-) (y) \right)^2 + 2 \left(v_\lambda (x) - v_\lambda (y) \right)\left( w^+(x) - w^+(y) \right). \notag
\end{align}
From \eqref{sec4-eq011} and \eqref{sec4-eq11}, we deduce that
\begin{align}\label{sec4-eq12}
    \mathcal{I}_\lambda (w) - \mathcal{I}_\lambda (0) &\geq \frac{1}{2} \|v_\lambda + w^- \|_{\mathcal{D}}^2 - \iint\limits_{\R^N \times \R^N} \frac{(v_\lambda (x) - v_\lambda (y)) (w^-(x) - w^-(y))}{|x-y|^{N+2s}}\,\dxy  - \lambda \int\limits_{\R^N} h(x)G(x,w)\,\dx - \frac{1}{2} \|v_\lambda \|_{\mathcal{D}}^2 \notag\\
    & = \frac{1}{2} \|v_\lambda + w^- \|_{\mathcal{D}}^2 - \lambda \int\limits_{\R^N} h(x)F(v_\lambda + w^-)\,\dx - \frac{1}{2} \|v_\lambda \|_{\mathcal{D}}^2 + \lambda \int\limits_{\R^N} h(x)F(v_\lambda)\,\dx \notag\\
    & \quad - \iint_{\R^N \times \R^N} \frac{(v_\lambda (x) - v_\lambda (y)) (w^-(x) - w^-(y))}{|x-y|^{N+2s}}\,\dxy +\lambda \int_{\R^N} h(x)f(v_\lambda) w^-(x)\,\dx \notag\\
    & = \mathcal{J}_\lambda (v_\lambda + w^-) - \mathcal{J}_\lambda (v_\lambda ) - \langle \mathcal{J}'_\lambda (v_\lambda ), w^- \rangle \notag\\
    &= \mathcal{J}_\lambda (v_\lambda + w^-) - \mathcal{J}_\lambda (v_\lambda ) \geq 0.
\end{align}
 The last inequality follows since $v_\lambda$ is a local minimizer for $\mathcal{J}_\lambda$ in $\Ds2$ (see Remark \ref{Rem-2}). Thus, $w\equiv 0$ is a local minimizer for $\mathcal{I}_\lambda$ in $\Ds2$.

Further, we claim that $-v_\lambda$ is a local minimizer for $\mathcal{I}_\lambda$. Notice that $G(x,- v_\lambda)= - F(v_\lambda) + f(v_\lambda)v_\lambda$, since $- v_\lambda \leq 0$. Recall that $v_\lambda$ is a critical point of $\mathcal{J}_\lambda$ and subsequently we get
\begin{align}\label{sec4-eq9}
    \mathcal{I}_\lambda (- v_\lambda) - \mathcal{I}_\lambda (0) &=  \|v_\lambda\|_{\mathcal{D}}^2 - \lambda \int_{\R^N} h(x)G(x,- v_\lambda)\,\dx - \frac{1}{2} \|v_\lambda\|_{\mathcal{D}}^2 \notag\\
    & = \|v_\lambda\|_{\mathcal{D}}^2 + \lambda \int_{\R^N} h(x)F(v_\lambda)\,\dx - \lambda \int_{\R^N} h(x)f(v_\lambda)v_\lambda\,\dx - \frac{1}{2} \|v_\lambda\|_{\mathcal{D}}^2 \notag \\
    &= \langle \mathcal{J}'_\lambda(v_\lambda), v_\lambda \rangle - \mathcal{J}_\lambda (v_\lambda) = - \mathcal{J}_\lambda (v_\lambda).
\end{align}
Taking $\delta_1>0$ small enough and $w \in B_{\delta_1}(-v_\lambda) = \left\{ v \in \Ds2 : \|v+ v_\lambda\|_{\mathcal{D}}<\delta_1 \right\}$. Using \eqref{sec4-eq11} and the fact that $v_\lambda$ is a solution to \eqref{main-problem}, we get
\begin{align}\label{sec4-eq10}
    \delta_1^2 &> \|w+ v_\lambda\|_{\mathcal{D}}^2 \notag\\
    &\geq \|w^- + v_\lambda\|_{\mathcal{D}}^2 + \|w^+\|_{\mathcal{D}}^2 + 2\iint_{\R^N \times \R^N} \frac{(v_\lambda (x) - v_\lambda (y)) (w^+(x) - w^+(y))}{|x-y|^{N+2s}}\,\dxy \quad  \notag\\
    &= \|w^- + v_\lambda\|_{\mathcal{D}}^2 + \|w^+\|_{\mathcal{D}}^2 + 2 \lambda \int_{\R^N} h(x)f(v_\lambda)w^+(x)\,\dx \quad  \notag\\
    & \geq \|w^- + v_\lambda\|_{\mathcal{D}}^2, \text{ for all } w \in B_{\delta_1}(-v_\lambda).
\end{align}
From \eqref{sec4-eq12}, \eqref{sec4-eq9} and \eqref{sec4-eq10}, we deduce that
\begin{align*}
    \mathcal{I}_\lambda (w) - \mathcal{I}_\lambda (-v_\lambda) = \mathcal{I}_\lambda (w) - \mathcal{I}_\lambda (0) + \mathcal{I}_\lambda (0) - \mathcal{I}_\lambda (-v_\lambda) \geq \mathcal{J}_\lambda(v_\lambda + w^-) \geq \mathcal{J}_\lambda(0) =0.
\end{align*}
Finally, we get $- v_\lambda$ is a local minimizer of the functional $\mathcal{I}_\lambda$ in $\Ds2$.
\end{proof}
\begin{lemma}\label{lemma-4.7}
Let {\rm \ref{f1}-\ref{f3}} and {\rm \ref{h}} hold and $v_\lambda$ be a weak solution to \eqref{main-problem} as given in Proposition \ref{Prop-4.7}.    Then, $\mathcal{I}_\lambda$ satisfies the Palais Smale condition in $\Ds2$ i.e., if there is a sequence $\{w_n\} \subset \Ds2$ such that $\{\mathcal{I}_\lambda(w_n)\}$ is bounded in $\R$ and $\mathcal{I}'_\lambda(w_n) \rightarrow 0$ in the dual space $\left(\Ds2\right)^*$, then $\{w_n\}$ admits a convergent subsequence in $\Ds2$.
\end{lemma}
\begin{proof}
Let $\{w_n\} \subset \Ds2$ be a Palais Smale sequence for $\mathcal{I}_\lambda$. Then
\begin{align*}
    &\{\mathcal{I}_\lambda(w_n)\} \text{ is bounded in  } \R,\\
    &\mathcal{I}'_\lambda(w_n) \rightarrow 0, \text{ as } n \rightarrow \infty \text{ in the dual space } (\Ds2)^* .
\end{align*}
First we show that the sequence $\{w_n\}$ is bounded in $\Ds2$. Let us define
\begin{align*}
    \mathcal{A}(u,v)= \iint_{\R^N \times \R^N} \frac{(u(x)-u(y)) (v(x)-v(y))}{|x-y|^{N+2s}}\,\dxy.
\end{align*}
Notice that $\mathcal{A}(u,v)= \mathcal{A}(u,v^+) + \mathcal{A}(u,v^-)$, where $v^+= \max\{v,0\}, v^-=\min\{v,0\}$. From \eqref{sec4-eq8} and \eqref{sec4-eq11}, we get 
\begin{align}\label{sec4-eq13}
 \mathcal{I}_\lambda (w_n) &= \frac{1}{2}\|v_\lambda + w_n\|_{\mathcal{D}}^2- \lambda \int_{\R^N} h(x) G(x,w_n)\,\dx- \mathcal{A}(v_\lambda,w_n) \notag\\
    &\geq \frac{1}{2}\|v_\lambda + w_n^-\|_{\mathcal{D}}^2 +\frac{1}{2}\|w_n^+\|_{\mathcal{D}}^2 + \mathcal{A}(v_\lambda,w_n^+) - \lambda \int\limits_{\R^N} h(x) F(v_\lambda + w_n^{-})\,\dx + \lambda \int\limits_{\R^N} h(x) F(v_\lambda)\,\dx \notag\\
    &\qquad + \lambda \int\limits_{\R^N} h(x) f(v_\lambda (x))w_n^{-}(x)\,\dx  - \mathcal{A}(v_\lambda,w_n) \notag\\
   &\geq \frac{1}{2}\|v_\lambda + w_n^-\|_{\mathcal{D}}^2  - \lambda \int\limits_{\R^N} h(x) F(v_\lambda + w_n^-)\,\dx +\frac{1}{2}\|w_n^+\|_{\mathcal{D}}^2 - \mathcal{A}(v_\lambda,w_n^-) + \lambda \int\limits_{\R^N} h(x) f(v_\lambda (x))w_n^-(x)\,\dx  \notag\\ 
    &= \mathcal{J}_\lambda(v_\lambda + w_n^-) +\frac{1}{2}\|w_n^+\|_{\mathcal{D}}^2 - \langle\mathcal{J}'_\lambda(v_\lambda), w_n^- \rangle \notag\\
    &= \mathcal{J}_\lambda(v_\lambda + w_n^-) +\frac{1}{2}\|w_n^+\|_{\mathcal{D}}^2. 
\end{align}
The last equality follows using the fact that $v_\lambda$ is a weak solution to \eqref{main-problem}, which implies $\langle\mathcal{J}'_\lambda(v_\lambda), \varphi \rangle=0$, for all $\varphi \in \Ds2$. From \eqref{Exis-1-eq-1}, we recall that
\begin{align*}
    \mathcal{J}_\lambda (u) \geq \frac{1}{4} \|u\|_{\mathcal{D}}^2 - \lambda C_\lambda \|h\|_1.
\end{align*}
We combine the above inequality with \eqref{sec4-eq13} to get 
\begin{align*}
 \mathcal{I}_\lambda (w_n) \geq \mathcal{J}_\lambda(v_\lambda + w_n^-) +\frac{1}{2}\|w_n^+\|_{\mathcal{D}}^2 \geq \frac{1}{4} \|v_\lambda + w_n^-\|_{\mathcal{D}}^2 +\frac{1}{2}\|w_n^+\|_{\mathcal{D}}^2 - \lambda C_\lambda \|h\|_1,
\end{align*}
which further implies 
\begin{align*}
 \frac{1}{4} \|v_\lambda + w_n^-\|_{\mathcal{D}}^2 +\frac{1}{2}\|w_n^+\|_{\mathcal{D}}^2 \leq \mathcal{I}_\lambda (w_n) + \lambda C_\lambda \|h\|_1. 
\end{align*}
The above inequality infers that $\{v_\lambda + w_n^-\}$ and $\{w_n^+\}$ are bounded in $\Ds2$. Consequently, $\{w_n\}$ is also bounded in $\Ds2$. By reflexivity of $\Ds2$, there exists $w \in \Ds2$ such that up to a subsequence $w_n \rightharpoonup w$ in $\Ds2$, $w_n \rightarrow w$ in $L_h^2(\R^N)$ and 
\begin{align}\label{sec4-eq14}
    \langle \mathcal{I}'_\lambda (w_n), w_n-w \rangle \rightarrow 0 \text{ as } n \rightarrow \infty.
\end{align}
Further we prove that $w_n \rightarrow w$ in $\Ds2$. We observe that 
\begin{align*}
     \langle \mathcal{I}_{\lambda}'(w_n), w_n-w \rangle
     &= \mathcal{A}(v_\lambda+w_n, w_n-w) - \lambda \int_{\R^N} h(x) g(x,w_n)(w_n-w)\,\dx -  \mathcal{A}(v_\lambda, w_n-w)\\
     &=\mathcal{A}(v_\lambda+w_n, w_n-w) - \lambda \int_{\{w_n<0\}} h(x) f(v_\lambda + w_n)(w_n-w)\,\dx\\
    &\qquad + \lambda \int_{\{w_n<0\}} h(x) f(v_\lambda)(w_n-w)\,\dx -  \mathcal{A}(v_\lambda, w_n-w)\\
\end{align*}
\begin{align*}
    &= \langle \mathcal{J}_{\lambda}'(v_\lambda + w_n), w_n-w \rangle - \langle \mathcal{J}_{\lambda}'(v_\lambda), w_n-w \rangle + \lambda \int_{\{w_n \geq 0\}} h(x) f(v_\lambda + w_n)(w_n-w)\,\dx\\
    & \qquad - \lambda \int_{\{w_n \geq 0\}} h(x) f(v_\lambda)(w_n-w)\,\dx.
\end{align*}
Using the compact embedding $\Ds2 \hookrightarrow L_{h}^q({\R^N}), q\in [1,2_s^*)$ (cf. Proposition \ref{embedding}-(ii)), we have
\begin{align}\label{sec4-eq15}
    \int_{\{w_n \geq 0\}} h(x) f(v_\lambda + w_n)(w_n-w)\,\dx \rightarrow 0 \text{ as } n \rightarrow \infty,
\end{align}
and 
\begin{align}\label{sec4-eq16}
    \int_{\{w_n \geq 0\}} h(x) f(v_\lambda)(w_n-w)\,\dx \rightarrow 0 \text{ as } n \rightarrow \infty.
\end{align}
Using \eqref{sec4-eq15} and \eqref{sec4-eq16}, we write
\begin{align}\label{sec4-eq17}
    \langle \mathcal{I}_{\lambda}'(w_n), w_n-w \rangle &= \langle \mathcal{J}_{\lambda}'(v_\lambda + w_n), w_n-w \rangle - \langle \mathcal{J}_{\lambda}'(v_\lambda), w_n-w \rangle + o_n(1) \notag\\
    &= \langle \mathcal{J}_{\lambda}'(v_\lambda + w_n), w_n-w \rangle + o_n(1).
\end{align}
From \eqref{sec4-eq14} and \eqref{sec4-eq17}, we obtain
\begin{align}\label{sec4-eq18}
    \langle \mathcal{J}_{\lambda}'(v_\lambda + w_n), w_n-w \rangle \rightarrow 0 \text{ as } n \rightarrow \infty.
\end{align}
We define $Q_\lambda : \Ds2 \rightarrow \R$ and $\psi_\lambda: \Ds2 \rightarrow \R$ as follows
\begin{align*}
    Q_\lambda (w) = \frac{1}{2} \|v_\lambda + w\|_{\mathcal{D}}^2,\quad \text{ and }\quad
    \psi_\lambda (w) = \lambda \int_{\R^N} h(x) F(v_\lambda +w)\,\dx.
\end{align*}
The functionals $Q_\lambda$ and $\psi_\lambda$ are Frech\'et differentiable on $\Ds2$ and we get
\begin{align}\label{sec4-eq19}
    \langle \psi_{\lambda}'(w_n), w_n-w \rangle = \lambda \int_{\R^N} h(x) f(v_\lambda + w_n)(w_n-w)\,\dx \rightarrow 0 \text{ as } n \rightarrow \infty.
\end{align}
Moreover, we also have
\begin{align*}
    \langle Q_{\lambda}'(w_n), w_n-w \rangle = \langle \mathcal{J}_{\lambda}'(v_\lambda + w_n), w_n-w \rangle + \langle \psi_{\lambda}'(w_n), w_n-w \rangle.
\end{align*}
From \eqref{sec4-eq18} and \eqref{sec4-eq19}, we obtain
\begin{align}\label{sec4-eq20}
    \langle Q_{\lambda}'(w_n), w_n-w \rangle \rightarrow 0 \text{ as } n \rightarrow \infty \quad\text{ i.e.,}\quad \mathcal{A}(v_\lambda+w_n, w_n-w) \rightarrow 0.
\end{align}
Since $v_\lambda \in \Ds2$ and $w_n \rightharpoonup w$, we also deduce from above that
\begin{align}\label{sec4-eq21}
  \mathcal{A}(w_n, w_n-w)=  \iint_{\R^N \times \R^N} \frac{( w_n(x) - w_n(y)) ((w_n-w))(x)-(w_n-w)(y))}{|x-y|^{N+2s}}\,\dxy \rightarrow 0 \text{ as } n \rightarrow \infty.
\end{align}
Thus, we write
\begin{align*}
    \|w_n - w\|_{\mathcal{D}}^2 &= \mathcal{A}(w_n, w_n-w)  - \mathcal{A}(w, w_n-w).
\end{align*}
The first term on the right hand side of the above equality tends to $0$ as $n \rightarrow \infty$ due to \eqref{sec4-eq21} and the second term goes to $0$ using the definition of weak convergence $w_n \rightharpoonup w$ in $\Ds2$. Finally, we obtain that $\|w_n - w\|_{\mathcal{D}} \rightarrow 0$ as $n \rightarrow \infty$. This completes the proof.
\end{proof}
\begin{definition}
    Given a ball $B\subset Y$, where $Y$ is a Banach space and $c \in \R$, the functional $I: Y \rightarrow \R$ is said to satisfies the Palais Smale condition around $B$ at level $c$ denoted by $(PS)_{B,c}$, if for every $\{u_n\} \subset Y$ such that
    \begin{itemize}
        \item $I(u_n) \rightarrow c$ in $\R$;
        \item $\text{dist}(u_n,B) \rightarrow 0$, where $\text{dist}(u_n,B)=\inf\limits_{v \in B}\|u_n-v\|_{Y}$;
        \item $I'(u_n) \rightarrow 0$ in the dual space $Y^*$,
    \end{itemize}
    the sequence $\{u_n\}$ admits a convergent subsequence in $Y$.
\end{definition}
In the following, we define the sub-level sets of $I$ and the set of all critical points of $I$ at level $c$.
\begin{align*}
    I^\alpha = \{x \in Y : I(x) \geq \alpha \} \text{ and } K_c(I) = \{x \in Y : I'(x)=0 \text{ and } I(x)=c\}.
\end{align*}
Further we state the Classical Linking Theorem due to \cite{Linking-11-Du}.
\begin{theorem}[\bf Classical Linking Theorem]\label{Linking-Theorem}
    Let $I: Y \rightarrow \R$ be a $C^1$ functional satisfying the $(PS)_{B,c}$ condition for all $c \in \R$, where $B= \{x \in Y: \|x-x_i\|_{Y}=r\} \subset I^\alpha$ with $x_i \in \{x_1,x_2\} \in Y$, $x_1 \neq x_2$, $r \in (0, \|x_2-x_1\|_Y)$ and $\alpha = \max\{I(x_1), I(x_2)\}$. Then
    \begin{align*}
        c_* = \inf_{\gamma \in \Gamma} \max_{s \in [0,1]}I(\gamma(s)) \geq \max\{I(x_1), I(x_2)\}, \text{ where } \Gamma = \{\gamma \in C([0,1],Y) : \gamma (0) = x_1, \gamma (1) = x_2 \},
    \end{align*}
    is a critical value for the functional $I$ and $K_{c_*}(I) \setminus \{x_1,x_2\} \neq \emptyset$.
\end{theorem}
\begin{lemma}\label{Lemma-4.10}
    There exists a critical point $w_\lambda \in \Ds2$ of $\mathcal{I}_\lambda$ such that $w_\lambda \leq 0$ in $\R^N$ and $w_\lambda \not\equiv 0,-v_\lambda$.
\end{lemma}
\begin{proof}
Without loss of generality, we assume that $\mathcal{I}_\lambda (-v_\lambda) \geq \mathcal{I}_\lambda (0)$. Since $-v_\lambda$ is a local minimizer to $\mathcal{I}_\lambda$ in $\Ds2$, there exists $\delta \in (0, \|v_\lambda\|_{\mathcal{D}})$ such that
\begin{align}\label{sec4-eq22}
    \mathcal{I}_\lambda (v) \geq \mathcal{I}_\lambda (-v_\lambda) := \alpha, \text{ for all } v \in \overline{B}_\delta(-v_\lambda).
\end{align}
Define $R = \{v \in \Ds2 : \|v+ v_\lambda\|_{\mathcal{D}}=\delta \}$ and $\mathcal{I}_\lambda^\alpha = \{v \in \Ds2 : \mathcal{I}_\lambda(v) \geq \alpha \}$. Notice that $R= \partial B_\delta (-v_\lambda) \subset \mathcal{I}_\lambda^\alpha$. Moreover, we have
\begin{align*}
    \inf\{\mathcal{I}_\lambda(v) : \|v+ v_\lambda\|_{\mathcal{D}}=\delta \} \geq \mathcal{I}_\lambda (-v_\lambda) = \max\{\mathcal{I}_\lambda(0), \mathcal{I}_\lambda (-v_\lambda)\}.
\end{align*}
By Lemma \ref{lemma-4.7} the functional  $\mathcal{I}_\lambda$ verify the $(PS)_{R,c}$-conditions for all $c \in \R$. Now we take
\begin{align*}
    c_\lambda = \inf_{\gamma \in \Gamma} \max_{s \in [0,1]}\mathcal{I}_\lambda(\gamma(s)) \geq \max\{\mathcal{I}_\lambda(0), \mathcal{I}_\lambda(-v_\lambda)\},  
\end{align*}
where $\Gamma = \{\gamma \in C([0,1],\Ds2) : \gamma (0) = 0, \gamma (1) = -v_\lambda \}$. Using the Classical Linking Theorem \ref{Linking-Theorem}, we have that $c_\lambda$ is a critical value to $\mathcal{I}_\lambda$. Thus, there exists $w_\lambda \in \Ds2$ which is different to $0$ and $-v_\lambda$ such that
\begin{align}\label{seq4-eq23}
    \mathcal{I}_\lambda (w_\lambda) = c_\lambda, \text{ and } \langle \mathcal{I}'_\lambda (w_\lambda), \varphi \rangle =0, \text{ for all } \varphi \in \Ds2.
\end{align}
Now we only need to show that $w_\lambda$ is non-positive in $\R^N$. We write $w_\lambda = w_\lambda^+ + w_\lambda^-$, where $w_\lambda^- = \min\{w_\lambda,0\}, w_\lambda^+ = \max\{w_\lambda,0\}$.  Notice that $w_\lambda^+ \geq 0$ and $w_\lambda^- \leq 0$. Putting $\varphi = w_\lambda^+$ in the above equation \eqref{seq4-eq23}. Then we have
\begin{align}\label{sec4-eq24}
  0=  \langle \mathcal{I}'_\lambda (w_\lambda), w_\lambda^+ \rangle = \iint_{\R^N \times \R^N} \frac{(w_\lambda(x)-w_\lambda(y)) (w_\lambda^+(x)-w_\lambda^+(y))}{|x-y|^{N+2s}}\,\dxy - \lambda \int_{\R^N} h(x) g(x,w_\lambda)w_\lambda^+(x)\,\dx.
\end{align}
Further, we have
\begin{align}\label{sec4-eq25}
    (w_\lambda(x)-w_\lambda(y)) (w_\lambda^+(x)-w_\lambda^+(y)) &= (w_\lambda^+(x)-w_\lambda^+(y))^2 + (w_\lambda^-(x)-w_\lambda^-(y)) (w_\lambda^+(x)-w_\lambda^+(y)) \notag\\
    & = (w_\lambda^+(x)-w_\lambda^+(y))^2 - w_\lambda^-(y)w_\lambda^+(x) - w_\lambda^-(x)w_\lambda^+(y) \notag\\
    & \geq (w_\lambda^+(x)-w_\lambda^+(y))^2.
\end{align}
From \eqref{sec4-eq24} and \eqref{sec4-eq25}, and using the fact that $\text{supp}(w_\lambda^+) \subset D:=\{x\in \R^N : w_\lambda(x)> 0\}$, we get the following estimate
\begin{align*}
    0 \geq  \iint_{\R^N \times \R^N} \frac{(w_\lambda^+(x)-w_\lambda^+(y))^2}{|x-y|^{N+2s}}\,\dxy - \lambda \int_{D} h(x) g(x,w_\lambda)w_\lambda(x)\,\dx = \iint_{\R^N \times \R^N} \frac{(w_\lambda^+(x)-w_\lambda^+(y))^2}{|x-y|^{N+2s}}\,\dxy.
\end{align*}
Thus, we get $\|w_\lambda^+\|_{\mathcal{D}}=0$, which implies that $w_\lambda^+ \equiv 0$ in $\R^N$ since $w_\lambda^+ \in \Ds2$. This completes the proof.
\end{proof}

\begin{proof}[Proof of Theorem \ref{second-solution}]
We prove that $\widetilde{v}_\lambda=v_\lambda + w_\lambda$ is the second solution to \eqref{main-problem} and is also positive. Recall that $w_\lambda \leq 0$ in $\R^N$ and we have
\begin{align*}
    g(x,w_\lambda(x)) = f(v_\lambda(x)+ w_\lambda(x)) - f(v_\lambda(x)) \text{ in } \R^N.
\end{align*}
Notice that $\langle \mathcal{I}'_\lambda (w_\lambda), \varphi \rangle=0$, for all $\varphi \in \Ds2$. Therefore, we have
\begin{align*}
&\iint_{\R^N \times \R^N} \frac{(v_\lambda(x) + w_\lambda(x) -v_\lambda(y) - w_\lambda(y)) (\varphi(x)-\varphi(y))}{|x-y|^{N+2s}}\,\dxy-  \iint_{\R^N \times \R^N} \frac{(v_\lambda(x)-v_\lambda(y)) (\varphi(x)-\varphi(y))}{|x-y|^{N+2s}}\,\dxy\\
&= \lambda \int_{\R^N} h(x) g(x,w_\lambda)\varphi\,\dx\\ &=\lambda \int_{\R^N} h(x)\left[f(v_\lambda + w_\lambda) - f(v_\lambda) \right] \varphi\,\dx . 
\end{align*}
Using that $v_\lambda$ is a solution to \eqref{main-problem}, the above identity reduces to
\begin{align*}
\iint_{\R^N \times \R^N} \frac{(v_\lambda(x) + w_\lambda(x) -v_\lambda(y) - w_\lambda(y)) (\varphi(x)-\varphi(y))}{|x-y|^{N+2s}}\,\dxy =\lambda \int_{\R^N} h(x)f(v_\lambda + w_\lambda) \varphi\,\dx, 
\end{align*}
for all $\varphi \in \Ds2$. We conclude that $\widetilde{v}_\lambda=v_\lambda + w_\lambda$ is also a solution to \eqref{main-problem} and $\widetilde{v}_\lambda  \not\equiv v_\lambda$. Moreover, the solution $\widetilde{v}_\lambda$ is positive using the strong maximum principle.
\end{proof}
\subsection{Examples of \texorpdfstring{$f$}{TEXT} and \texorpdfstring{$h$}{TEXT}} In this subsection, we provide nontrivial examples of $f$ and $h$ which satisfies the required assumption in our analysis.
\begin{example}\label{Example on f}
    Let us consider a function $f: \R \rightarrow \R$ as follows:  
\begin{align*}
\displaystyle f(t) =   \begin{cases}
        \frac{t^p}{t^{p-\sigma} + 1}, &\quad \text{for } t > 0,\\
        0 , &\quad \text{for } t \leq 0,
    \end{cases}
\end{align*}
for some parameters \( p > 1 \) and \( \sigma \in (0,1) \). Then $f$ satisfies all the assumptions \ref{f1}--\ref{f5}.
\end{example}
\begin{example}\label{Example on h}
    Let $\eta\in C_c^\infty([0,\infty))$ be a cutoff with $0\le\eta\le1$ and
\[
\eta(r)=1\ \text{ for }0\le r\le1,\quad \eta(r)=0\ \text{ for }r\ge2. 
\]
Define
\[
P(r):=\eta(r)+(1-\eta(r))(1+r)^{-(N+2s)},\quad r\ge0.
\]
Then $P\in C([0,\infty),[0,\infty))$, %$P(r)\sim (1+r)^{-(N+2s+1)}$ as $r\to\infty$, and hence 
and $P\in L^\infty([0,\infty))\cap L^1([0,\infty))$. We define
$h(x):=P(|x|),\ x\in  \mathbb{R}^N.$ Clearly $h\in L^\infty(\mathbb{R}^N)\cap L^1(\mathbb{R}^N)$ and $0<h(x)\le P(|x|)$ for all $x$. Set
\[
I(x):=\int_{\mathbb{R}^N}P(|y|)\,|x-y|^{2s-N}\,\dy .
\]
We verify \ref{P1} and \ref{P2}. %Fix constants $C_k>0$ whose values may change from line to line.
Notice that $P(|y|)\lesssim |y|^{-(N+2s)}$ for $y \neq 0$. For any fixed $R>>1$ large enough, we split the integral as
\[
I(x)=\int_{|x-y|\le R}P(|y|)\,|x-y|^{2s-N}\,\dy
+\int_{|x-y|>R}P(|y|)\,|x-y|^{2s-N}\,\dy.
\]
For $|x|\leq 1$, we estimate the first integral as
\begin{align*}
    \int_{|x-y|\le R}P(|y|)\,|x-y|^{2s-N}\,\dy \leq C_1\|P\|_\infty \int_{0}^{R}r^{2s-N} r^{N-1}\,\dr = \frac{C\|P\|_\infty R^{2s}}{2s}.  
\end{align*}
The second integral is estimated as follows:
\begin{align*}
    \int_{|x-y|>R}P(|y|)\,|x-y|^{2s-N}\,\dy &=  \int_{|z|>R}P(|x-z|)\,|z|^{2s-N}\,\dz\\
    &\leq C_2 \int_{|z|>R}\frac{|z|^{2s-N}}{|x-z|^{N+2s}}\,\dz\\
    & \leq C_2 \int_{|z|>R}\frac{|z|^{2s-N}}{\left(|z|-|x|\right)^{N+2s}}\,\dz\\
    &\leq C_3 \int_{R}^{\infty} r^{2s-N} r^{N-1} r^{-N-2s}\,\dr = C_4 R^{-N}.
\end{align*}
% The first term is finite for $|x|\leq R$ because the kernel $|x-y|^{2s-N}$ is locally integrable. For the second term use $P(|y|)\lesssim |y|^{-(N+2s+1)}$ for large $|y|$ and $|x-y|\gtrsim|y|$ to obtain an absolutely convergent tail. 
The last step follows since $|z|>R$, where $R$ is very large. Thus $\sup_{|x|\le 1}I(x)<\infty$ and so $I\in L^\infty_{\mathrm{loc}}(\mathbb{R}^N)$.
For $|x|>1$, we write
\begin{align}\label{e-0}
    \int_{\mathbb{R}^N}P(|y|)\,|x-y|^{2s-N}\,\dy = \int_{|x-y| > |x|/2 }P(|y|)\,|x-y|^{2s-N}\,\dy + \int_{|x-y| \leq |x|/2}P(|y|)\,|x-y|^{2s-N}\,\dy.
\end{align}
The first integrals is estimated as follows:
\begin{align}\label{e-1}
    \int_{|x-y| > |x|/2 }P(|y|)\,|x-y|^{2s-N}\,\dy &\leq \left(\frac{1}{2}\right)^{2s-N} |x|^{2s-N} \int_{|x-y| > |x|/2 }P(|y|)\,\dy \notag\\
    & \leq C_5 |x|^{2s-N}.
\end{align}
Now we estimate the second integral. For $|x-y| \leq \frac{|x|}{2}$, we set $z=x-y$. Then $|x-z| \geq |x|-|z| \geq \frac{|x|}{2}$. Thus we have
\begin{align}\label{e-2}
    \int_{|x-y| \leq |x|/2 }P(|y|)\,|x-y|^{2s-N}\,\dy &= \int_{|z| \leq |x|/2 }P(|x-z|)\,|z|^{2s-N}\,\dz \notag\\
    & \leq C_6 \int_{|z| \leq |x|/2 } \frac{|z|^{2s-N}}{|x-z|^{N+2s}}\,\dz \notag \\
    &\leq C_7 \left(\frac{2}{|x|}\right)^{N+2s} \int_{0}^{\frac{|x|}{2}} r^{2s-N} r^{N-1}\,\dr = C_8 |x|^{-N} \leq C_8 |x|^{2s-N}.
\end{align}
From \eqref{e-0},\eqref{e-1} and \eqref{e-2}, we deduce that
 \[
 I(x)\le C_9|x|^{2s-N},\quad\text{for }|x|>1,
 \]
which is \ref{P2}. Together with the local boundedness this implies $\sup_{x\in\mathbb{R}^N}I(x)<\infty$, i.e. \ref{P1} holds.
\end{example}

\begin{example}\label{New Example}
For some $\delta>0$, we consider the function
    \begin{align}\label{e-3}
        h(x) = \frac{1}{(1+|x|)^{N} \, (\log(e+|x|))^{1+\delta}},\, x\in \R^N \text{ and }
P(r) = \frac{1}{(1+r)^{N} (\log(e+r))^{1+\delta}},\, r \geq 0.
\end{align}
Then $h$ is a positive continuous function on $\R^N$ and it satisfies $0<h(x) \leq P(|x|)$, for all $x \in \R^N$. Notice that $h \in L^\infty(\R^N)$. Moreover, we have
\begin{align*}
    \int_{\mathbb{R}^N} h(x)\,dx 
&= |\mathbb{S}^{N-1}|\int_0^\infty 
\frac{r^{N-1}}{(1+r)^N (\log(e+r))^{1+\delta}}\,\dr\\
&=|\mathbb{S}^{N-1}|\int_0^R 
\frac{r^{N-1}}{(1+r)^N (\log(e+r))^{1+\delta}}\,\dr + |\mathbb{S}^{N-1}|\int_R^\infty 
\frac{r^{N-1}}{(1+r)^N (\log(e+r))^{1+\delta}}\,\dr\\
& \leq C_1 \int_0^R r^{N-1}\,\dr + C_2 \int_R^\infty 
\frac{1}{r (\log(r))^{1+\delta}}\,\dr,
\end{align*}
where $R$ is some fixed large number. Both the integrals on the right hand side are finite which infers $h \in L^1(\R^N)$. 

For $|x|\leq 1$, we split the integral $I$ as follows:
\begin{align*}
    I(x) = \int_{\mathbb{R}^N} P(|y|)\,|x-y|^{2s-N}\,\dy = \int_{|y|\leq 2|x|} P(|y|)\,|x-y|^{2s-N}\,\dy +  \int_{|y|>2|x|} P(|y|)\,|x-y|^{2s-N}\,\dy.
\end{align*}
The first integral on the right hand side is finite since the function $P$ is bounded and the kernel $|x-y|^{2s-N}$ is locally integrable. The second integral is estimated as follows:
\begin{align*}
    \int_{|y|>2|x|} P(|y|)\,|x-y|^{2s-N}\,\dy
\leq C_2 \int_{|y|>2|x|} P(|y|)\,|y|^{2s-N}\,\dy
\leq C_3 \int_{2|x|}^\infty r^{2s-N-1}(\log r)^{-1-\delta}\,\dr.
\end{align*}
The last integral converges uniformly in $x$. Hence, the integral $I$ is finite and $\sup_{|x|\leq 1}I(x)< + \infty$.

For $|x|> 1$, we split the integral $I$ as follows:
\begin{align}\label{e-4}
   I(x)= \int_{|x-y| > |x|/2 }P(|y|)\,|x-y|^{2s-N}\,\dy + \int_{|x-y| \leq |x|/2}P(|y|)\,|x-y|^{2s-N}\,\dy.
\end{align}
The first integrals is estimated as follows:
\begin{align}\label{e-5}
    \int_{|x-y| > |x|/2 }P(|y|)\,|x-y|^{2s-N}\,\dy &\leq \left(\frac{1}{2}\right)^{2s-N} |x|^{2s-N} \int_{|x-y| > |x|/2 }P(|y|)\,\dy \notag \\
    & \leq C_4 |x|^{2s-N}.
\end{align}
Now we estimate the second integral. For $|x-y| \leq \frac{|x|}{2}$, we set $z=x-y$. Then $|x-z| \geq |x|-|z| \geq \frac{|x|}{2}$. Thus we have
\begin{align}\label{e-6}
    \int_{|x-y| \leq |x|/2 }P(|y|)\,|x-y|^{2s-N}\,\dy &= \int_{|z| \leq |x|/2 }P(|x-z|)\,|z|^{2s-N}\,\dz \notag\\
    & = \int_{|z| \leq |x|/2 } \frac{|z|^{2s-N}}{(1+|x-z|)^{N} \, (\log(e+|x-z|))^{1+\delta}}\,\dz \notag\\
    &\leq \left(\frac{2}{|x|} \right)^N \int_{0}^{|x|/2} r^{2s-N} r^{N-1}\,\dr \notag\\
    &= C_5 |x|^{2s-N}.
\end{align}
From \eqref{e-4}, \eqref{e-5} and \eqref{e-6}, we get
\begin{align}\label{e-7}
    I(x)\le C_6|x|^{2s-N}, \qquad |x|> 1,
\end{align}
which proves \ref{P2}. Moreover, \ref{P1} holds by combining the boundedness of $I$ for $|x|\leq 1$ and \eqref{e-7} for $|x|>1$.

For every $\gamma>0$ we have
\begin{align*}
    \frac{h(r)}{(1+r)^{-(N+\gamma)}} 
= (1+r)^{\gamma} (\log(e+r))^{-1-\delta}
\rightarrow \infty
\quad \text{as } r\to\infty.
\end{align*}
Thus, there does not exist any positive constant $C_7$ such that
\begin{align*}
    h(x) \le C_7(1+|x|)^{-(N+\gamma)}, \text{ for all } |x|>1.
\end{align*}
Hence, the assumptions \ref{P1} and \ref{P2} on $h$ are weaker than the assumption of Carl, Perera, and Tehrani~\cite[Theorem~5.2]{Carl-Perera-Tehrani-2025-Arxiv} (see also \eqref{eq-h}).
\end{example}

\section*{Acknowledgment}
J. Abrantes dos Santos was partially supported by CNPq, Brazil, grant 304774/2022-7, and FAPESQ, Brazil, grant 3031/2021. R. Kumar acknowledges TIFR-Centre for Applicable Mathematics, India for the Institute Postdoctoral Fellowship to conduct this research. A. Sarkar expresses his heartfelt gratitude for the DST-INSPIRE Grant DST/INSPIRE/04/2018/002208, funded by the Government of India.

%%%%%%%%%%%%%%%%%%%%%%%%%%%%%%%%%%%%%%%%%%%%%%%%%%%%%%%%%%%%
%\bibliography{ref}
%\bibliographystyle{abbrv}

\end{document}